\documentclass[a4paper,twosided,12pt]{amsart}

\usepackage[latin1]{inputenc}
\usepackage{hyperref}
\usepackage{amsxtra}
\usepackage{amssymb}
\usepackage{amsmath,amsthm}
\usepackage[all]{xypic}
\usepackage{epsfig}
\xyoption{dvips}

\theoremstyle{definition}
\newtheorem{ntn}{Notation}[subsection]
\newtheorem{dfn}[ntn]{Definition}
\newtheorem{rem}[ntn]{Remark}
\newtheorem{exa}[ntn]{Example}
\theoremstyle{plain}

\newtheorem{lem}[ntn]{Lemma}

\newtheorem{prp}[ntn]{Proposition}
\newtheorem{thm}[ntn]{Theorem}
\newtheorem{cor}[ntn]{Corollary}

\theoremstyle{remark}

\DeclareMathAlphabet{\mathds}{U}{dsrom}{m}{n}
\DeclareMathAlphabet{\mathsc}{U}{rsfs}{m}{n}

\DeclareMathOperator{\Aut}{Aut}
\DeclareMathOperator{\cha}{char}
\DeclareMathOperator{\id}{id}
\DeclareMathOperator{\Sch}{Sch}

\DeclareMathOperator{\Spec}{Spec}
\DeclareMathOperator{\Spf}{Spf}

\DeclareMathOperator{\Hm}{Hom}

\DeclareMathOperator{\Is}{Isom}
\DeclareMathOperator{\Pic}{Pic}

\DeclareMathOperator{\NS}{NS}
\DeclareMathOperator{\DefF}{Def}
\DeclareMathOperator{\CSpin}{CSpin}
\DeclareMathOperator{\CSp}{CSp}

\DeclareMathOperator{\SO}{SO}
\DeclareMathOperator{\GL}{GL}
\DeclareMathOperator{\PrG}{PGL}

\newcommand{\Q}{\mathbb{Q}}
\newcommand{\N}{\mathbb{N}}
\newcommand{\Z}{\mathbb{Z}}

\renewcommand{\P}{\mathbb{P}}
\newcommand{\R}{\mathbb{R}}
\newcommand{\A}{\mathbb{A}}
\newcommand{\C}{\mathbb{C}}
\newcommand{\Gm}{\mathbb{G}}

\newcommand{\Kg}{\mathbb{K}}

\renewcommand{\O}{\mathcal O}
\renewcommand{\L}{\mathcal L}

\newcommand{\M}{\mathcal M}

\newcommand{\X}{\mathcal X}

\renewcommand{\a}{\alpha}
\renewcommand{\l}{\lambda}

\newcommand{\G}{\Gamma}

\newcommand{\e}{\epsilon}

\newcommand{\Fk}{\mathcal F}
\newcommand{\Av}{\mathcal A}

\newcommand{\lr}{\rightarrow}

\title{Moduli Stacks of Polarized K3 Surfaces in Mixed Characteristic}
\author{Jordan Rizov}
\address{Mathematisch Instituut\\ P.O. Box 80.010\\ 3508 TA Utrecht\\ The Netherlands}
\email{rizov@math.uu.nl\footnote{Author's current address is: ABN AMRO N.V., Enterprise Risk Modelling, Gustav Mahlerlaan 10, 1082 PP Amsterdam, The Netheralds\\{\it E-mail address:} jordan.rizov@nl.abnamro.com}}
\keywords{K3 surfaces, Moduli spaces}
\subjclass[2000]{14J28, 14D22}

\begin{document}
\maketitle
\begin{abstract}
In this note we define moduli stacks of (primitively) polarized K3 spaces. We show that they are representable by Deligne-Mumford stacks over $\Spec(\Z)$. Further, we look at K3 spaces with a level structure. Our main result is that the moduli functors of K3 spaces with a primitive polarization of degree $2d$ and a level structure are representable by smooth algebraic spaces over open parts of $\Spec(\Z)$. To do this we use ideas of Grothendieck, Deligne, Mumford, Artin and others.

These results are the starting point for the theory of complex multiplication for K3 surfaces and the definition of Kuga-Satake abelian varieties in positive characteristic given in our Ph.D. thesis \cite{JR-Thesis}.
\end{abstract}
\section*{Introduction}
In this note we will consider moduli spaces of K3 surfaces with a polarization. For a natural number $d$ and an algebraically closed field $k$, a K3 surface with a polarization of degree $2d$ over $k$ is a
pair $(X,\L)$ consisting of a K3 surface $X$ over $k$ and an ample line bundle
$\L$ on $X$ with self intersection number $(\L,\L) = 2d$. The moduli space of polarized K3 surfaces with certain level structure over $\C$ is constructed as an open subspace of the Shimura variety associated with $\SO(2,19)$. Over $\Z$ we use techniques developed by Artin to show the existence of such spaces.

In various places in the literature one finds detailed accounts on coarse moduli schemes of primitively polarized complex K3 surfaces. We outline in Section \ref{ModStackSection} two approaches to the theory, one via geometric invariant theory (\cite{Vie-M}) and another via periods of complex K3 surfaces (\cite[Expos\'e XIII]{Ast-K3} and \cite[\S 1]{Fri-Torelli}). Here we take up a different point of view and work with moduli stacks rather than with coarse moduli schemes. In this way, our exposition is closer to \cite{Olsson-K3} where moduli stacks of primitively polarized K3 surfaces and their compactifictions over $\Q$ are constructed. We define the categories $\Fk_{2d}$ and $\M_{2d}$ of primitively polarized (respectively
polarized) K3 surfaces of degree $2d$ over $\Z$ and show that they are Deligne-Mumford stacks over $\Z$.

For various technical reasons we will need to work with algebraic spaces rather than with Deligne-Mumford stacks. In the case of abelian varieties one
introduces level $n$-structures using Tate modules and considers moduli functors of polarized abelian varieties with level
$n$-structure for $n \in \N, \ n\geq 3$. These functors are representable by schemes. We adopt a similar strategy in order to define moduli functors which are
representable by algebraic spaces. For a certain class of
compact open subgroups $\Kg$ of $\SO(2,19)(\A_f)$ we introduce the notion of a level $\Kg$-structure on K3 surfaces using their second \'etale cohomology groups. Further, we introduce moduli spaces $\Fk_{2d,\Kg}$ of
primitively polarized K3 surfaces with level $\Kg$-structure and show that these are smooth algebraic spaces over $\Spec(\Z[1/N_\Kg])$ where $N_\Kg \in \N$
depends on $\Kg$. These moduli spaces are finite unramified covers of $\Fk_{2d}$. Important examples of level structures are spin level $n$-structures. These are level structures defined by the images of
some principal level $n$-subgroups of $\CSpin(2,19)(\A_f)$ under the adjoint representation homomorphism $\CSpin(2,19) \lr \SO(2,19)$. We denote the corresponding moduli
space by $\Fk_{2d,n^{\rm sp}}$.

Let us outline briefly the contents of this note. In the first few sections we review some basic properties of K3 surfaces. Then we continue with the study of the representability of Picard and automorphism functors arising from K3 surfaces. The core of the problems discussed here is Section \ref{ModStackSection} in which we define various moduli functors of polarized K3 surfaces and prove that those define Deligne-Mumford stacks. In Section \ref{LevelStructures} we define level structures on K3 surfaces associated to compact open subgroups of $\SO(2,19)(\A_f)$. In the last section we show that the moduli functors of primitively polarized K3 surfaces with level structure are representable by algebraic spaces.
\newline
\newline
{\bf Notations}
\newline
\newline
We write $\hat \Z$ for the profinite completion of $\Z$. We denote by $\A$ the ring of ad\`eles of $\Q$ and by $\A_f = \hat \Z \otimes \Q$ the ring of finite ad\`eles of $\Q$. Similarly, for a number field $E$ we denote by $\A_E$ and $\A_{E,f}$ the ring of ad\`eles and the ring of finite ad\`eles of $E$.

If $A$ is a ring, $A \lr B$ a ring homomorphism then for any $A$-module ($A$-algebra etc.) $V$ we will denote by $V_B$ the $B$-module ($B$-algebra etc.) $V\otimes_A B$.

For a ring $A$ we denote by $({\rm Sch}/A)$ the category of schemes over $A$. We will write ${\rm Sch}$ for the category of schemes over $\Z$.

By a variety over a field $k$ we will mean a separated, geometrically integral scheme of finite type over $k$. For a variety $X$ over $\C$ we will denote by $X^{\rm an}$ the associated analytic variety. For an algebraic stack $\mathcal F$ over a scheme $S$ and a morphism of schemes $S' \lr S$ we will denote by $\mathcal F_{S'}$ the product $\mathcal F\times_S S'$ and consider it as an algebraic stack over $S'$.

A superscript ${}^0$ indicates a connected component for the Zariski topology. For an algebraic group $G$ will denote by $G^0$ the connected component of the identity. We will use the superscript ${}^+$ to denote connected components for other topologies.

Let $V$ be a vector space over $\Q$ and let $G \hookrightarrow \GL(V)$ be an algebraic group over $\Q$. Suppose given a full lattice $L$ in $V$ (i.e., $L\otimes \Q = V$). Then $G(\Z)$ and $G(\hat \Z)$ will denote the abstract groups consisting of the elements in $G(\Q)$ and $G(\A_f)$ preserving the lattices $L$ and $L_{\hat \Z}$ respectively.
\newline
\newline
{\bf Acknowledgments}
\newline
\newline
This note contains the results of Chapter 1 of my Ph.D. thesis \cite{JR-Thesis}. I thank my advisors, Ben Moonen and Frans Oort for their help, their support and for everything I have learned from them. I would like to thank Bas Edixhoven and Gerard van der Geer for pointing out some mistakes and for their valuable suggestions. I thank the Dutch Organization for Research N.W.O. for the financial support with which my thesis was done.
\section{Basic Results}
\subsection{Definitions and Examples}
We will briefly recall some basic notions concerning families of K3
surfaces.
\begin{dfn}
Let $k$ be a field. A non-singular, proper surface $X$ over $k$ is called a \emph{K3 surface} if $\Omega^2_{X/k} \cong \O_X$ and $H^1(X,\O_X) = 0$.
\end{dfn}
Note that a K3 surface is automatically projective. Let us give some basic examples one can keep in mind:
\begin{exa}\label{DoubCov}
Let $S$ be a non-singular sextic curve in $\P^2_k$ where $k$ is a field and consider a
double cover i.e., a finite generically \'etale morphism, $\pi \colon X \lr \P^2_k$ which is ramified along $S$. Then $X$ is a
K3 surface.
\end{exa}
\begin{exa}\label{ComplInt}
{\it Complete intersections}: Let $X$ be a smooth surface which is a
complete intersection of $n$ hypersurfaces of degree $d_1,\dots,d_n$ in $\P^{n+
2}$ over a field $k$. The adjunction formula shows that $\Omega^2_{X/k} \cong \O_X(d_1 +
\dots + d_n - n - 3)$. So a necessary condition for $X$ to be a K3 surface is
$d_1 + \dots + d_n = n + 3$. The first three possibilities are:
$$
\begin{matrix}
  n = 1 & d_1 = 4 \\
  n = 2 & d_1 = 2, d_2 = 3 \\
  n = 3 & d_1 = d_2 = d_3 = 2.
\end{matrix}
$$
For a complete intersection $M$ of dimension $n$ one has that
$H^i(M,\O_M(m)) = 0$ for all $m \in \Z$ and $1 \leq i \leq n-1$. Hence in those
three cases we have $H^1(X, \O_X) = 0$ and therefore $X$ is a K3 surface.
\end{exa}
\begin{exa}\label{KummerS}
Let $A$ be an abelian surface over a field $k$ of characteristic different from
2. Let $A[2]$ be the kernel of the multiplication by-2-map,
let $\pi \colon \tilde A \lr A$ be the blow-up of $A[2]$ and let $\tilde E$
be the exceptional divisor. The automorphism
$[-1]_A$ lifts to an involution $[-1]_{\tilde A}$ on $\tilde A$. Let $X$ be the quotient
variety of $\tilde A$ by the group of automorphisms $\{\id_{\tilde
A}, [-1]_{\tilde A}\}$ and denote by $\iota: \tilde A \lr X$ the
quotient morphism. It is a finite map of degree 2. We have the following diagram
$$
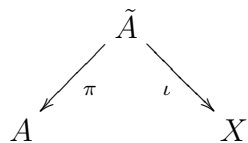
\xymatrix{
& \tilde A \ar[dl]^\pi \ar[dr]_\iota & \\
A & & X
}
$$
of morphisms over $k$. The variety $X$ is a K3 surface and it is called the \emph{Kummer surface} associated to $A$.
\end{exa}
\begin{dfn}\label{DefK3Sp}
By a \emph{K3 scheme} over a base scheme $S$ we will mean a scheme $X$ and a proper and smooth morphism $\pi \colon X \lr S$ whose geometric fibers are K3 surfaces. A \emph{K3 space over a scheme $S$} is an algebraic space $X$ together with a proper and smooth morphism $\pi \colon X \lr S$ such that there is an \'etale cover $S' \lr S$ of $S$ for which $\pi' \colon X' = X\times_S S' \lr S'$ is a K3 scheme.
\end{dfn}
If $\pi \colon X \lr S$ is a K3 space, then $\pi_* \O_X = \O_S$.
Indeed, this is true since $\pi$ is proper and its geometric fibers are reduced and connected.
\begin{rem}
A K3 space $X$ over $S$ is usually defined as an algebraic space $X$ together with a proper and smooth morphism $\pi \colon X \lr S$ such that for every geometric point $s \in S$ the fiber $X_s$ is a K3 surface. In this note we will restrict ourselves to Definition \ref{DefK3Sp} above. The reason is that for this class of K3 spaces one can easily see that certain automorphism functors of K3 spaces are representable by schemes (cf. Theorem \ref{AutGrK3}). We do not know if this holds in general.
\end{rem}
\subsection{Ample Line Bundles on K3 Surfaces}
In order to construct the moduli stacks of polarized K3 spaces one
needs a number of results on ample line bundles. We give them below.
\begin{dfn}
Let $X$ be a K3 surface over a field $k$. The self-intersection index $(\L,\L)_X$ of a line bundle $\L$ on $X$ will be called its \emph{degree}. A line bundle $\L$ on $X$ is called \emph{primitive} if $\L\otimes \bar k$ is is not a positive power of a line bundle on $X_{\bar k}$.
\end{dfn}
\begin{thm}\label{albprop} 
Let $X$ be a K3 surface over a field $k$.
\begin{enumerate}
\item [(a)] If $\L$ is a line bundle on $X$, then $(\L,\L)$ is even. If $\L$ is ample and $d := (\L,\L)/2$, then the Hilbert polynomial of $\L$ is given by $h_\L(t) =dt^2 + 2$.
\item [(b)] Suppose $\L$ is an ample bundle. Then $\L$ is effective and $H^i(X,\L) = 0$ for $i>0$. Further, $\L^n$ is generated by global sections if $n \geq 2$ and is very ample if $n \geq 3$.
\end{enumerate}
\end{thm}
\begin{proof}
(a) First note that, by Serre duality, $h^2(\O_X) = h^0(\Omega^2_{X/k}) = h^0(\O_X) = 1$. Since $h^0(\O_X) = 1$ we find that $\chi(\O_X) = 2$. Hirzebruch-Riemann-Roch gives
\begin{displaymath}
\begin{split}
\chi(\L) &= \chi(\O_X) + \frac{1}{2} \cdot \bigl((\L,\L) - (\L,\Omega^2_{X/k})\bigr) \\
	 &= 2 + \frac{1}{2}\cdot (\L,\L) \\
\end{split}
\end{displaymath}
as $\Omega^2_{X/k}$ is trivial. Hence $(\L,\L) = 2d$ is even. If $\L$ is ample then its Hilbert polynomial is $h_{\L}(t) =dt^2 + 2$

(b) By Serre duality and the fact that $\Omega^2_{X/k} \cong \O_X$ we have $h^i(\L) = h^{2-i}(\L^{-1})$. In particular $h^2(\L) = h^0(\L^{-1}) = 0$ as an anti-ample bundle is not effective. Since $d := (\L,\L)/2 > 0$ it follows that $h^0(\L) = d + 2 + h^1(\L) > 0$, so $\L$ is effective. For the remaining assertions we refer to \cite{Saint-D}, Section 8.
\end{proof}
\begin{exa}
Let $\pi \colon X \lr \P^2$ be a double cover of $\P^2$ as in Example \ref{DoubCov}. The line bundle $\L = \pi^*\O_{\P^2}(1)$ is ample and one has that $(\L,\L)_X = 2 (\O_{\P^2}(1),\O_{\P^2}(1))_{\P^2} = 2$. Hence any K3 surface $X$ which is a double cover of $\P^2$ ramified along a non-singular sextic curve has an ample line bundle $\L$ of degree 2. 
\end{exa}
\begin{exa}
Let $X \subset \P^{n+2}$ be a K3 surface which is obtained as a complete intersection of multiple degree $(d_1,d_2,\dots,d_n)$; see Example \ref{ComplInt}. Then $\O_{X}(1)$ degree $d_1d_2\cdots d_n$. Note that the equality $d_1 + d_2 + \cdots + d_n = n+3$ implies that at least one of the $d_i$ is even.
\end{exa}
Note that if $\pi \colon X \lr S$ is a K3 scheme over a connected base $S$ then for a line bundle $\L$ on $X$ the intersection index $(\L_{\bar s},\L_{\bar s})_{X_{\bar s}}$ is constant for any $\bar s$. This follows from the fact that $\pi$ is flat and the relation $(\L_{\bar s}, \L_{\bar s})_{X_{\bar s}} = 2\chi(\L_{\bar s}) - 4$.
\begin{lem}\label{RelAmpl}
Let $\pi \colon X \lr S$ be a K3 scheme and let $\L$ be a line
bundle on $X$ which is fiberwise ample on $X$ i.e., $\L_{\bar s}$ is ample on $X_{\bar s}$
for every geometric point $\bar s \in S$. Let $2d = (\L_{\bar s},\L_{\bar s})_{X_{\bar s}}$ for any point $\bar s \in S$. Then $\pi_*\L^n$ is a locally free sheaf of of rank $dn^2+2$ and $\L^n$ is relatively very ample over $S$ if $n \geq 3$.
\end{lem}
\begin{proof}
By Theorem \ref{albprop} (b) we have that for all $\bar s \in S$ the group $H^1(X_{\bar s},\L_{\bar s}^n)$ is trivial. It follows from \cite[Ch. III, \S 7]{EGA}, that $\pi_*\L^n$ is a locally free sheaf and that $\pi_*\L_{\bar s} \cong H^0(\X_{\bar s}, \L^n_{\bar s})$. The rank statement follows from Theorem \ref{albprop} (a). By part (a) of Theorem \ref{albprop} one sees that for every geometric point $\bar s \in S$ and any $n \geq 3$ the line bundle $\L^n_{\bar s}$ gives a closed immersion $X_{\bar s} \hookrightarrow \P(\pi_*\L^n_{\bar s})$ over $\kappa(\bar s)$. Hence the morphism $X \hookrightarrow \P(\pi_*\L^n)$ induced by $\L^n$ is a closed immersion. This finishes the proof. 
\end{proof}
\section{Cohomology Groups of K3 Surfaces}
\subsection{Quadratic Lattices Related to Cohomology Groups of K3 Surfaces}\label{QLatK3}
In this section we introduce some notations which will be used in the sequel. Let $U$ be the hyperbolic plane and denote by $E_8$ the \emph{positive} quadratic lattice associated to the Dynkin diagram of type $E_8$ (cf. \cite[Ch. V, 1.4 Examples]{S-CA}). 
\begin{ntn}
Denote by $(L_0,\psi)$ the quadratic lattice $U^{\oplus 3} \oplus E_8^{\oplus 2}$. Further, let $(V_0,\psi_0)$ be the quadratic space $(L_0,\psi) \otimes_\Z \Q$. 
\end{ntn}
\noindent
We have that $L_0$ is a free $\Z$-module of rank 22. The form $\psi_\R$ has signature $(19+,3-)$ on $L_0 \otimes \R$.

Let $\{e_1,f_1\}$ be a basis of the first copy of $U$ in $L_0$ such that
$$
 \psi(e_1,e_1) = \psi(f_1,f_1) = 0\ \ \text{and}\ \ \psi(e_1,f_1) = 1.
$$
For a positive integer $d$ we consider the vector $e_1-df_1$ of $L_0$. It is a primitive vector i.e., the module $L_0/\langle e_1-df_1\rangle$ is free and we have that $\psi(e_1-df_1,e_1-df_1) = -2d$. The orthogonal complement of $e_1-df_1$ in $L_0$ with respect to $\psi$ is $\langle e_1+df_1\rangle \oplus U^{\oplus 2} \oplus E^{\oplus 2}_8$.
\begin{ntn}
Denote the quadratic sublattice $\langle e_1+df_1\rangle \oplus U^{\oplus 2} \oplus E^{\oplus 2}_8$ of $L_0$ by $(L_{2d},\psi_{2d})$. Further, we denote by $(V_{2d},\psi_{2d})$ the quadratic space $(L_{2d},\psi_{2d})\otimes_\Z \Q$. 
\end{ntn}
The signature of the form $\psi_{2d,\R}$ is $(19+,2-)$. We have that $\langle e_1-df_1\rangle \oplus L_{2d}$ is a sublattice of $L_0$ of index $2d$. The inclusion of lattices $i \colon L_{2d} \hookrightarrow L_0$ defines injective homomorphisms of groups
\begin{equation}\label{InjOHom}
i^{\rm ad} \colon \{g \in {\rm O}(V_0)(\Z)\ |\ g(e_1-df_1) = e_1-df_1\}\hookrightarrow {\rm O}(V_{2d})(\Z)
\end{equation}
and
\begin{equation}\label{InjSOHom}
i^{\rm ad} \colon \{g \in \SO(V_0)(\Z)\ |\ g(e_1-df_1) = e_1-df_1\} \hookrightarrow \SO(V_{2d})(\Z).
\end{equation}
Let $L_{2d}^*$ denote the dual lattice $\Hm(L_{2d},\Z)$. Then the bilinear form $\psi_{2d}$ defines an embedding $L_{2d} \hookrightarrow L_{2d}^*$ and we denote by $A_{2d}$ the factor group $L_{2d}^*/L_{2d}$. It is an abelian group of order $2d$ (\cite[\S 2, Lemma]{L-P}). One can extend the bilinear form $\psi_{2d}$ on $L_{2d}$ to a $\Q$-valued form on $L_{2d}^*$ and define
\begin{displaymath}
q_{2d} \colon A_{2d} \lr \Q/2\Z
\end{displaymath}
defined by
\begin{displaymath}
 q_{2d}(x + L_{2d}) = \psi_{2d}(x,x) + 2\Z
\end{displaymath} 
for any $x \in L_{2d}^*$. Let ${\rm O}(q_{2d})$ denote the group of isomorphisms of $A_{2d}$ preserving the form $q_{2d}$. Then one has a natural homomorphism $\tau \colon {\rm O}(V_{2d})(\Z) \lr {\rm O}(q_{2d})$. It is shown in \cite{Nik-QLat} that 
\begin{displaymath}
 i^{\rm ad}\bigl(\{g \in {\rm O}(V_0)(\Z)\ |\ g(e_1-df_1) = e_1-df_1\}\bigr) = \ker(\tau).
\end{displaymath} 
\subsection{De Rham Cohomology}
Let $X$ be a K3 surface over a field $k$. The following proposition will play an essential role when studying deformations of K3 surfaces (Section \ref{DefTh}). We will use it also to show that the automorphism group $\Aut(X)$ of a K3 surface is reduced (see Theorem \ref{AutGrK3} below).
\begin{prp}\label{HodgeNum}
 If $X$ is a K3 surface over a field $k$, then
\begin{enumerate}
\item [(a)] The Hodge-de Rham spectral sequence
 $$
  E_1^{i,j} = H^j(X,\Omega_{X/k}^i) \Longrightarrow H^{i+j}_{DR}(X,k)
 $$
 degenerates at $E_1$. For the Hodge numbers $h^{i,j} = \dim_k H^j(X,\Omega_{X/k}^i)$ of $X$ we have 
\begin{gather*}
   h^{1,0} = h^{0,1} = h^{2,1} = h^{1,2} = 0 \\
   h^{0,0} = h^{2,0} = h^{0,2} = h^{2,2} = 1 \\
   h^{1,1} = 20. \\
\end{gather*}
\item [(b)] Let $\Theta_{X/k} = \Omega_{X/k}^{1 \vee}$ be the tangent bundle of $X$. Then
 $H^i(X,\Theta_{X/k}) = 0$ for $i = 0$ and $2$ and $\dim_k H^1(X,\Theta_{X/k}) = 20$.
\end{enumerate}
\end{prp}
\begin{proof}
If $k$ has characteristic zero, then one may assume that $k = \C$ and the proposition follows from \cite[\S 1, Prop. 1.2]{L-P}. The case $\cha(k) = p >0$ is treated in \cite[Prop. 1.1]{Del-K3}.
\end{proof}
\begin{rem} 
Part (b) of the proposition is classical in the case $k = \C$. The proof in the general case is due to Rudakov and Shafarevich. It can be reformulated in following way: There exist no non-trivial regular vector fields on a K3 surface (cf. \cite[\S 6, Thm. 7]{RSh-IMorK3}).
\end{rem}
\subsection{Betti Cohomology}\label{BettiCohom}
Let $X$ be a complex K3 surface. Then the Betti cohomology groups
$H_B^i(X,\Z)$ are free $\Z$-modules of rank $1,0,22,0,1$ for
$i=0,1,2,3,4$ respectively. One has a non-degenerate bilinear form (given by the Poincar\'e duality pairing):
$$
 \psi \colon H_B^2(X,\Z)(1) \times H_B^2(X,\Z)(1) \lr \Z
$$
given by
$$
 \psi(x,y) = - {\rm tr} (x \cup y)
$$
where $x \cup y$ is the cup product of $x$ and $y$ and ${\rm tr} \colon H^4_B(X,\Z(2)) \lr \Z$ is the trace map. It has signature $(19+,3-)$ over $\R$. The quadratic lattice
$\bigl(H_B^2(X,\Z)(1), \psi\bigr)$ is isometric to $(L_0,\psi)$ (cf. Section \ref{QLatK3}). For proofs of those results we refer to \cite[\S1, Prop. 1.2]{L-P}.

The group $H^2_B(X,\Z)$ carries a natural $\Z$-Hodge structure (which we will abbreviate as $\Z$-HS) of type $\{(2,0), (1,1),(0,2)\}$ with $h^{2,0} = h^{0,2} =
1$ and $h^{1,1} = 20$ as we see from Proposition \ref{HodgeNum}.

For a complex K3 surface $H^1(X,\O_X)$ is trivial so the first Chern class map
$$
 c_1 \colon \Pic(X) \lr H^2_B(X,\Z)(1)
$$ 
is injective. Exactly in the same way we
see that for a K3 space $\pi \colon X \lr S$, where $S$ is a scheme over $\C$, one has a short exact
sequence of sheaves
$$
 0 \lr R^1\pi_*^{\rm an}\O_X^* \lr R^2\pi_*^{\rm an}\Z(1)
$$
as $R^1\pi_*^{\rm an}\O_X$ is trivial.

\begin{ntn}
Let $\L$ be an ample line bundle on $X$. We denote by $P_B^2(X,\Z)(1)$ the orthogonal complement of $c_1(\L)$ with respect to $\psi$. It
is a free $\Z$-module of rank 21 called the \emph{primitive part} (or the \emph{primitive cohomology group}) of $H_B^2(X,\Z)(1)$ with
respect to $c_1(\L)$. The restriction of $\psi$ defines a non-degenerate bilinear form:
\begin{displaymath}
 \psi_\L \colon P^2_B(X,\Z)(1) \times P^2_B(X.\Z)(1) \lr \Z.
\end{displaymath}
\end{ntn}
The group $P_B^2(X,\Z(1))$ carries a natural $\Z$-HS induced by the one on $H^2_B(X,\Z(1))$
of type $\{(-1,1), (0,0) , (1,-1)\}$ with $h^{-1,1} = h^{1,-1} = 1$ for which $\psi_\L$ is a
polarization.
\begin{rem}\label{PrimLB}
Let $\L$ be an ample line bundle for which $(\L,\L)_X=2d$ and assume that it is primitive. 
Let $\{e_1,f_1\}$ be a basis of the first copy of $U$ in $L_0$ as in Section \ref{QLatK3}. By \cite[Exp. IX, \S1, Prop. 1]{Ast-K3} one can find an isometry 
\begin{displaymath}
a \colon \bigl(H^2_B(X,\Z(1)),\psi\bigr) \lr L_0
\end{displaymath}
such that $a(c_1(\L)) = e_1 - df_1$. Therefore $a$ induces an isometry 
\begin{displaymath}
a \colon \bigl(P^2_B(X,\Z(1)), \psi_\L\bigr) \lr (L_{2d},\psi_{2d}).
\end{displaymath}
\end{rem}
\subsection{\'Etale Cohomology}
Let $k$ be a field of characteristic $p \geq 0$ and fix a prime $l$ which
is different from $p$. Suppose given a K3 surface $X$ over $k$. Then the \'etale cohomology group
$H_{\rm et}^i(X_{\bar k},\Z_l)$ is a free $\Z_l$-module of rank $1,0,22,0,1$
for $i = 0,1,2,3,4$. One sees this in the following way: If $k$ has characteristic
zero, then the claim follows from the corresponding result for Betti
cohomology and the comparison theorem between Betti and \'etale cohomology
(\cite[Ch. III, \S 3, Thm. 3.12]{Mil-EC}). Assume that $p > 0$. By \cite[\S 1, Cor. 1.8]{Del-K3} there
exists a discrete valuation ring $R$ with residue field $\bar k$ and a smooth lift $\X$ over $R$ of $X$.
If $\eta$ is the generic point of $\Spec(R)$, then by the smooth base change theorem for \'etale cohomology (\cite[Ch. VI, \S 4, Cor. 4.2]{Mil-EC}) one has that 
\begin{equation}\label{etcohrel}
 H^i_{\rm et}(X_{\bar k},\Z/l^n\Z) \cong H^i_{\rm et}(\X_{\bar \eta}, \Z/l^n\Z)
\end{equation}
for every $i=0,\dots, 4$ and every $n$. Hence $H^i_{\rm et}(X_{\bar k},\Z_l) \cong H^i_{\rm et}(\X_{\bar \eta}, \Z_l)$ and we deduce the claim from the characteristic zero result.

Further, one has a non-degenerate bilinear form
\begin{displaymath}
 \psi_{\Z_l} \colon H_{\rm et}^2(X_{\bar k},\Z_l)(1) \times H_{\rm et}^2(X_{\bar k},\Z_l)(1) \lr \Z_l
\end{displaymath}
given by
$$
 \psi_{\Z_l} (x,y) = -\text{tr}_{\Z_l}(x\cup y)
$$
where $\text{tr}_{\Z_l} \colon H^4_{\rm et}(X_{\bar k},\Z_l)(2) \lr \Z_l$ is the trace
isomorphism. This is simply Poincar\'e duality for \'etale cohomology (\cite[Ch.
VI, \S 11, Cor. 11.2]{Mil-EC}).

The Kummer short exact sequence of \'etale sheaves on $X$
$$
 1 \lr \boldsymbol{\mu}_{l^n} \lr \mathbb G_m \lr \mathbb G_m \lr 1
$$
gives an exact sequence of cohomology groups 
\begin{displaymath}
H^1_{\rm et}(X_{\bar k}, \boldsymbol{\mu}_{l^n}) \lr H_{\rm et}^1(X_{\bar k}, \mathbb G_m) \lr H_{\rm et}^1(X_{\bar k}, \mathbb G_m) \lr H^2_{\rm et}(X_{\bar k}, \boldsymbol{\mu}_{l^n}).
\end{displaymath}
By \eqref{etcohrel} the group $H^1_{\rm et}(X_{\bar k},\boldsymbol{\mu}_{l^n})$ is trivial we have an injection 
\begin{displaymath}
0 \lr \Pic(X)/l^n \Pic(X) \lr  H^2_{\rm et}(X_{\bar k}, \boldsymbol{\mu}_{l^n}).
\end{displaymath}
Taking the projective limit over $n$ one sees that the first Chern class map 
$$
 c_1 \colon  \Pic(X) \otimes_\Z \Z_l \hookrightarrow H_{\rm et}^2(X_{\bar k}, \Z_l)(1)
$$ 
is injective. In particular, since
$H^2_{\rm et}(X_{\bar k},\Z_l(1))$ is free, $\Pic(X)$ has no $l$-torsion for any $l$ different from $p$.

Similarly, if $\pi \colon X \lr S$ is a K3 space then one can consider the long exact
sequence of higher direct images, coming from the Kummer sequence
\begin{displaymath}
 R^1_{\rm et}\pi_*\boldsymbol{\mu}_{l^n} \lr R^1_{\rm et}\pi_*\mathbb G_m \lr R^1_{\rm et}\pi_*\mathbb
 G_m \lr R^2_{\rm et}\pi_*\boldsymbol{\mu}_{l^n}.
\end{displaymath}
Further, since the stalk of $R^1_{\rm et}\pi_*\boldsymbol{\mu}_{l^n}$ at any geometric point of $S$
is zero (one uses here the proper base change theorem), the sheaf itself is zero (\cite[Ch. II, \S 2, Prop. 2.10]{Mil-EC}). Hence passing again to the projective limit over
$n$ we obtain the exact sequence of $\Z_l$-sheaves 
\begin{displaymath}
 0 \lr R^1_{\rm et}\pi_*\mathbb G_m \otimes \Z_l \lr R^2_{\rm et}\pi_*\Z_l(1).
\end{displaymath}
\begin{ntn}
Let $\L$ be a primitive ample line bundle on $X$ with $(\L,\L)_X = 2d$. Denote by $P^2_{\rm et}(X_{\bar k}, \Z_l(1))$ the \emph{primitive part} of $H^2_{\rm et}(X_{\bar k},\Z_l)(1)$ with respect to
$c_1(\L)$ i.e., the orthogonal complement of $c_1(\L)$
in $H^2_{\rm et}(X_{\bar k},\Z_l)(1)$ with respect to $\psi_{\Z_l}$. Denote the restriction of $\psi_{\Z_l}$ to $P^2_{\rm et}(X_{\bar k},\Z_l(1))$ by $\psi_{\L,\Z_l}$.
\end{ntn}
If $k$ has characteristic 0, then by the comparison theorem between Betti and \'etale cohomology one
has that $\bigl(H^2_{\rm et}(X_{\bar k},\Z_l(1)),\psi_{\Z_l}\bigr)$ is isometric to $\bigl(H^2_B(X_\C,\Z(1)),\psi\bigr)\otimes_\Z \Z_l$ which is isometric to $(L_0,\psi) \otimes_\Z \Z_l$. Moreover since the comparison isomorphism respects algebraic cycles, the same holds for
the primitive parts with respect to $\L$ i.e., we have that $\bigl(P^2_{\rm et}(X_{\bar k},\Z_l(1)),
\psi_{\L,\Z_l}\bigr) \cong (L_{2d}, \psi_{2d}) \otimes_\Z \Z_l$.

Assume that $\cha(k) = p>0$. Then the pair $(X,\L)\otimes \bar k$ has a lift
$(\X, \L)$ over a discrete valuation ring $R$ with $\cha(R) =0$ and with residue field $\bar k$
(see \cite[\S 1, Cor. 1.8]{Del-K3}). Using the same argument
as above one concludes that 
\begin{displaymath}
H^i_{\rm et}(X_{\bar k}, \Z_l)(m) \cong H^i_{\rm et}(\X_{\bar \eta}, \Z_l)(m)
\end{displaymath}
and that $\bigl(H^2_{\rm et}(X_{\bar k}, \Z_l)(m), \psi_{\Z_l}\bigr)$
is isometric to $\bigl(H^2_{\rm et}(\X_{\bar \eta}, \Z_l)(m), \psi_{\Z_l}\bigl)$, where $\eta$ is the
generic point of $\Spec(R)$. Consequently the two quadratic lattices
$\bigl(P^2_{\rm et}(X_{\bar k},\Z_l(1)), \psi_{\L,\Z_l}\bigr)$ and $\bigr(P^2_{\rm et}(\X_{\bar \eta},\Z_l(1)), \psi_{\L_{\bar \eta}}\bigr)\otimes_\Z \Z_l$ are also isometric. Thus, if $\L$ is primitive, then there is an isometry 
\begin{displaymath}
a \colon \bigl(H^2_{\rm et}(X_{\bar k},\Z_l(1)), \psi_l\bigr) \lr L_0 \otimes \Z_l
\end{displaymath}
such that $a(c_1(\L)) = e_1-df_1$. It induces an isometry 
\begin{displaymath}
a \colon \bigl(P^2_{\rm et}(X_{\bar k}, \Z_l(1)), \psi_{\L,\Z_l}\bigr) \lr (L_{2d},\psi_{2d})\otimes \Z_l.
\end{displaymath}
\begin{rem}
Let $k$ be a  field of characteristic $p$. We make the following notations 
\begin{displaymath}
\hat \Z^{(p)} := \prod_{l\ne p} \Z_l\ \ \ \text{and}\ \ \ \A_f^{(p)} = \hat \Z^{(p)} \otimes \Q.
\end{displaymath} 
In the sequel we will be considering \'etale cohomology with $\hat \Z^{(p)}$ or $\A_f^{(p)}$ coefficients. Then we have that for a K3 surface over a field $k$ one has isometries
$$
 \bigl(H^2_{\rm et}(X_{\bar k},\hat \Z^{(p)}(1)), \psi_f\bigr) \cong (L_0,\psi)\otimes_\Z \hat \Z^{(p)}
$$
and for a primitive ample line bundle $\L$ of degree $2d$ on $X$ one has
$$
 \bigl(P^2_{\rm et}(X_{\bar k},\hat \Z^{(p)}(1)), \psi_{\L,f}\bigr) \cong (L_{2d},\psi_{2d})\otimes_\Z \hat \Z^{(p)}.
$$ 
Here $\psi_f$ and $\psi_{\L,f}$ are the corresponding bilinear forms coming from the Poincar\'e duality on $H^2_{\rm et}(X_{\bar k}, \hat \Z^{(p)}(1))$.
\end{rem}
\subsection{Crystalline Cohomology}
Let $k$ be a perfect field of characteristic $p>0$ and let $W =
W(k)$ be the ring of Witt vectors with coefficients in $k$. Consider a K3 surface
$X$ over $k$. Then by \cite[Prop. 1.1]{Del-K3} the crystalline
cohomology group $H^i_{\rm cris}(X/W)$ is a free $W$-module of rank
1, 0, 22, 0, 1 for $i =$ 0, 1, 2, 3, 4 respectively.
We consider next the crystalline Chern class map
$$
 c_1 \colon \Pic(X) \lr H^2_{\rm cris}(X/W).
$$ 
As pointed out in \cite[Appendice, Rem. 3.5]{Del-K3} the Chern class map defines
an injection 
$$
 c_1 \colon \NS(X_{\bar k}) \otimes_\Z \Z_p \hookrightarrow H^2_{\rm cris}(X/W(\bar k))
$$ 
where $\NS(X_{\bar k}) = \Pic(X_{\bar k}) / \Pic^0(X_{\bar k})$ is the N\'eron-Severi group of $X_{\bar k}$. In particular this means that the N\'eron-Severi group of $X_{\bar k}$  has no $p$-torsion.

If $K$ is the fraction field of $W$ then we shall denote by $H^i_{\rm cris}(X/K)$ the $K$-vector space $H^i_{\rm cris}(X/W)\otimes_{W} K$.
\section{Picard Schemes and Automorphisms of K3 Surfaces}
\subsection{Picard and N\'eron-Severi Groups of K3 Surfaces.}
In this section we will study Picard functors of K3 spaces. Those functors will play an important role in two aspects in the construction of moduli spaces of (primitively) polarized K3 surfaces. First, we will define (quasi-) polarizations on K3 surfaces using Picard spaces (cf. Definition \ref{PolDefinition} below). Later, in Section \ref{HilbertSchSection}, we will use Picard spaces in the construction of the Hilbert scheme parameterizing K3 subschemes of $\P^N$.

For a separated algebraic space $X$ over a scheme $S$ we denote by $\Pic(X)$ the group of isomorphism classes of invertible sheaves on $X$. Let $\pi \colon X \lr S$ be a K3 space and consider the \emph{relative Picard functor}
\begin{displaymath}
 \Pic_{X/S} \colon (\Sch/S)^0 \lr {\rm Groups}.
\end{displaymath} 
By definition it is the fppf-sheafification of the functor
$$
 P_{X/S} \colon (\Sch/S)^0 \lr {\rm Groups} \ \ \ \text{given\ by}\ \ \ T \mapsto \Pic(X \times_S T).
$$
For every $g \colon T \lr S$ we have that $\Pic_{X/S}(T) = H^0(T,R^1\pi'_*\Gm_m)$ where $\pi' \colon X \times_S T \lr T$ is the product morphism and all derived functors are taken with respect to the fppf-topology.
\begin{thm} 
For a K3 space $\pi \colon X \lr S$ the relative Picard functor $\Pic_{X/S}$ is represented by a separated algebraic space locally of finite presentation over $S$. 
\end{thm}
\begin{proof}
The representability follows form \cite[\S 7, Thm. 7.3]{Art-AFM1}. The proof of the separatedness property goes exactly in the same way as the proof of Theorem 3 in \cite[Ch. 8, \S 8.4]{NM}.
\end{proof}

Let $S = \Spec(k)$ be a spectrum of a field. Then $\Pic_{X/k}$ is represented by a group scheme (cf. \cite{O-SchDeP} or Lemma \ref{PicField} below) and shall denote by $\Pic^0_{X/k}$ its identity component. We set further
$$
 \Pic^{\tau}_{X/k} = \bigcup_{n > 0} n^{-1}\bigl(\Pic^0_{X/k}\bigr)
$$
where $n \colon \Pic_{X/k} \lr \Pic_{X/k}$ is the multiplication by $n$.
\begin{lem}\label{PicField}
Let $X$ be a K3 surface over a field $k$. Then $\Pic_{X/k}$ is represented by a
separated, smooth, zero dimensional scheme over $k$. In
particular $\Pic^0_{X/k}$ is trivial. Further, we have also that $\Pic^{\tau}_{X/k}$ is trivial.
\end{lem}
\begin{proof}
Combining Theorem 3 and Theorem 1, with $S = \Spec(k)$, of \cite[Ch. 8, \S
8.2]{NM} one concludes that $\Pic_{X/k}$ is representable by a separated scheme,
locally of finite type over $k$.

By Theorem 1 of \cite[Ch. 8, \S 8.4]{NM} one has that
$$
 \dim_k \Pic_{X/k} \leq \dim_k H^1(X,\O_X) = 0
$$
and as the equality holds in this case, $\Pic_{X/k}$ is smooth over $k$. This shows the validity of all assertions except for the claim about $\Pic^{\tau}_{X/k}$.

The scheme $\Pic^{\tau}_{X/k}$ is proper and of finite type over $k$ (cf. \cite[Ch. 8, Thm. 4]{NM}). Since its
dimension is zero it is a finite commutative group scheme over $k$. The
injectivity of the \'etale Chern class map shows that $\Pic(X)$ has no
$l$-torsion for $l \ne p$. By the first part of the lemma we have
that $\NS(X) = \Pic(X)$. Then the injectivity of the crystalline Chern class map shows
that $\Pic(X)$ has no $p$-torsion either. Thus $\Pic(X)$ is torsion free and
therefore $\Pic^{\tau}_{X/k}(\bar k)$ is trivial. Since in this case $\Pic^{\tau}_{X/k}$ is reduced we conclude it is trivial.
\end{proof}
If $X$ is a K3 surface over a field $k$, then $\NS(X) = \Pic(X)$, which follows from the fact that in this case $\Pic^0(X)$ is trivial. Hence $\Pic(X)$ is a free abelian group of rank at most $22$ (use \cite[Ch. V, \S 3, Cor 3.28]{Mil-EC}). If the characteristic of the ground field is zero, then  ${\rm rk}_\Z\Pic(X) \leq 20$.

Let $\pi \colon X \lr S$ be a K3 scheme. Define $\Pic^0_{X/S}$ and $\Pic^\tau_{X/S}$ as the subfunctors of $\Pic_{X/S}$ consisting of all elements whose restrictions to all fibers $X_s$ belong to $\Pic^0_{X_s/\kappa(s)}$ and $\Pic^\tau_{X_s/\kappa(s)}$ respectively.  
\begin{prp}\label{PicK3Spaces}
For a K3 scheme $\pi \colon X \lr S$ over a quasi-compact base $S$ one has that $\Pic_{X/S}$ is an algebraic space which is unramified over $S$. Further, we have that $\Pic^0_{X/S}$ and $\Pic^\tau_{X/S}$ are trivial.
\end{prp}
\begin{proof}
The first part of the proposition follows from the preceding lemma as it is enough to check the $\Pic_{X/S}$ is unramified in the case $S$ is a spectrum of a field. To prove the second part we notice that according to \cite[Ch. 8, \S 8.3, Thm. 4]{NM} we have open immersions $\Pic^0_{X/S} \hookrightarrow \Pic_{X/S}$ and $\Pic^\tau_{X/S} \hookrightarrow \Pic_{X/S}$. By Lemma \ref{PicField} above for every geometric point $\bar s \in S$ the subspaces $\Pic^0_{X_{\bar s}/\kappa(\bar s)}$ and $\Pic^\tau_{X_{\bar s}/\kappa(\bar s)}$ are trivial hence $\Pic^0_{X/S}$ and $\Pic^\tau_{X/S}$ are trivial.
\end{proof}
\begin{rem}\label{ModFuntConstr}
Let $\pi \colon X \lr S$ be a K3 scheme and let $\L$ and $\M$ be two line bundles on $X$. If $\L^n = \M^n$ for some $n \in \N$, then $\L$ is isomorphic to $\M \otimes \pi^*\mathcal N$ where $\mathcal N$ is a line bundle on $S$. Indeed, we have that $cl(\L)^n = cl(\M)^n$ in $\Pic_{X/S}$. Since $\Pic^\tau_{X/S}$ is trivial we have that the multiplication by $n$-morphism $[n]\colon \Pic_{X/S} \lr \Pic_{X/S}$ is an injective homomorphism of group schemes. Since $cl(\L)^n = cl(\M)^n$ we conclude that $cl(\L\otimes \M^{-1})$ is trivial, so $\M$ and $\L$ differ by an invertible sheaf coming from the base $S$ (\cite[Ch. 8, \S 8.1, Prop. 4]{NM}).
\end{rem}
\begin{rem}
It is easy to see that the statement of Proposition \ref{PicK3Spaces} remains true for K3 spaces. 
\end{rem}
A morphism of schemes $\pi \colon X \lr S$ is called \emph{strongly projective} (respectively \emph{strongly quasi-projective}) if there exists a locally free sheaf $\mathcal E$ on $S$ of constant finite rank such that $X$ is $S$-isomorphic to a closed subscheme (respectively a subscheme) of $\P(\mathcal E)$. 
\begin{lem}\label{PicDiv}
Let $S$ be a noetherian scheme and suppose given a K3 scheme $\pi \colon X \lr S$. If $\pi$ is a strongly projective morphism, then we have that
\begin{enumerate}
\item [(i)] for any $n \in \N$ the multiplication by $n$-morphism
$$
 [n] \colon \Pic_{X/S} \lr \Pic_{X/S}
$$
is a closed immersion of group schemes over $S$.
\item [(ii)] for any $\l \in \Pic_{X/S}(S)$ the set of points
$$
 S^o = \{ s \in S\ |\ \l_s\ {\rm is\ primitive\ on}\ X_s \}
$$
is open in $S$. 
\end{enumerate}
\end{lem}
\begin{proof}
{\it (i):} By definition we have a closed immersion $X \hookrightarrow \P(\mathcal E)$ for some locally free sheaf $\mathcal E$ on $S$. Let $\O_{X}(1)$ denote the pull-back of the canonical bundle $\O(1)$ on $\P(\mathcal E)$ via this inclusion. For a polynomial $\Phi \in \Q[t]$ let $\Pic^\Phi_{X/S}$ be the subfunctor of $\Pic_{X/S}$ which is induced by the line bundles $\L$ on $X$ with a given Hilbert polynomial $\Phi$ (with respect to $\O_X(1)$) on the fibers of $X$ over $S$. Then $\Pic_{X/S}^{\Phi}$ is representable by a strongly quasi-projective scheme over $S$ and $\Pic_{X/S}$ is the disjoint union of the open and closed subschemes $\Pic^\Phi_{X/S}$ for all $\Phi \in \Q[t]$. For a proof of this result we refer to \cite[Ch. 8, \S 8.2, Thm. 5]{NM}.

Since all schemes $\Pic^\Phi_{X/S}$ are quasi-compact we have that for a given $\Phi$ the image $[n](\Pic^\Phi_{X/S})$ is contained in a finite union $\bigcup_{i \in {\mathcal C}^n_\Phi} \Pic^{\Phi_i}_{X/S}$. We will show first that for a given $\Phi \in \Q[t]$ the morphism 
$$
 [n] \colon \Pic^\Phi_{X/S} \lr \bigcup_{i \in {\mathcal C}^n_\Phi} \Pic^{\Phi_i}_{X/S}
$$
is proper. As all schemes involved are noetherian we can apply the valuative criterion for properness. We may assume that $S$ is a spectrum of a discrete valuation ring $R$ and that $X$ admits a section over $S$ and let $\eta$ and $s$ be the generic and the special point of $S$. Under those assumptions any element of $\Pic_{X/S}$ comes from a class of a line bundle (\cite[Ch. 8, \S 8.1, Prop. 4]{NM}). To show that the restriction of $[n]$ to $\Pic_{X/S}^\Phi$ is proper we have to show that if $\L$ is a line bundle over the generic fiber $X_\eta$ of $X$, then $\L^n$ extends uniquely to a line bundle on $X$ which is a $n$-th power of a line bundle. This follows from \cite[Ch. 8, \S 8.4, Thm. 3]{NM} as both $\L$ and $\L^n$ extend uniquely over $X$.

Further, the morphism $[n] \colon \Pic_{X/S} \lr \Pic_{X/S}$ is an immersion of the corresponding topological spaces and as it is proper on every open and closed $\Pic^\Phi_{X/S}$, the image $[n](\Pic_{X/S})$ is closed in $\Pic_{X/S}$. We are left to show that the natural homomorphism of sheaves $\O_{\Pic_{X/S}} \lr [n]_*\O_{\Pic_{X/S}}$ is surjective. As this can be checked on stalks we see further that it is enough to show the surjectivity assuming that $S$ is a spectrum of a field. But under this condition the claim follows from Lemma \ref{PicField}. Indeed, $\Pic_{X/k}$ is a reduced, zero dimensional scheme. Hence all subschemes $\Pic^\Phi_{X/k}$ being reduced, quasi-projective and zero dimensional, are finite unions of points. Then the restrictions $[n] \colon \Pic^\Phi_{X/k} \lr \bigcup_{i \in {\mathcal C}^n_\Phi} \Pic^{\Phi_i}_{X/k}$ are closed immersions and hence $[n] \colon \Pic_{X/k} \lr \Pic_{X/k}$ is also a closed immersion. Therefore $\O_{\Pic_{X/k}} \lr [n]_*\O_{\Pic_{X/k}}$ is surjective.

{\it (ii):} We may assume that $S$ is connected. Then the intersection index $(\l_{\bar s},\l_{\bar s})$ is constant on $S$, say $(\l_{\bar s},\l_{\bar s}) = 2d$. For any natural number $n$ consider the closed subscheme $S_n$ of $S$ defined by the following Cartesian diagram
\begin{displaymath}
\xymatrix{
S_n \ar[d] \ar[r] & S \ar[d]^{\l} \\
\Pic_{X/S} \ar[r]^{[n]} & \Pic_{X/S}.
}
\end{displaymath}
Then the subset $S^o$ of $S$ can be identified with $S \setminus \bigcup_n S_n$ where the union is taken over all $n \in \N$ such that $n^2$ divides $d$. So it has a structure of an open subscheme of $S$. 
\end{proof}
\begin{rem}\label{EtPic}
Note that if $\pi \colon X \lr S$ is a K3 scheme, then the Picard functor $\Pic_{X/S}$ can be constructed using the \'etale topology on $S$ instead of the fppf-topology. In other words $\Pic_{X/S}$ is also the \'etale sheafification of $P_{X/S}$. This follows from the fact that $\pi$ is a proper morphism, using the Leray spectral sequence for $\pi$ and the sheaf $\Gm_m$. For a proof we refer to the comments on p. 203 in \cite[Ch. 8, \S 8.1]{NM}.
\end{rem}
\begin{exa}\label{NSKumS}
Let $A$ be an abelian surface over an algebraically closed field $k$ of characteristic different from
2 and let $X$ be the associated Kummer surface. Then one has that
$$
 \Pic(X)_\Q  = \NS(X)_\Q \cong \NS(A)^{[-1]_A}_\Q \oplus \Q^{\oplus
 16} 
$$ 
where $\NS(A)^{[-1]_A}$ denotes the elements of $\NS(A)$ invariant under the action of
$[-1]_A$. We refer to \cite[\S 3, Prop. 3.1]{Shi-SSK3} for a proof.
\end{exa}
\subsection{Polarizations of K3 Surfaces}\label{Polarizations}
Here we will define the notion of a polarization on a K3 space.
\begin{dfn}
Let $k$ be a field. A \emph{polarization} on a K3 surface $X/k$ is a global section $\l \in \Pic_{X/k}(k)$ which over $\bar k$ is the class of an ample line bundle $\L_{\bar k}$. The degree of $\L_{\bar k}$ is called the \emph{polarization degree} of $\l$. A \emph{quasi-polarization} on $X$ is a global section $\l \in \Pic_{X/k}(k)$ which over $\bar k$ comes from a line bundle $\L_{\bar k}$ with the following property: 
\begin{enumerate}
\item [(i)] $\L_{\bar k}$ is nef i.e., $\bigl(\L_{\bar k},\O_{X_{\bar k}}(C)\bigr) \geq 0$ for all irreducible curves in $X_{\bar k}$,
\item [(ii)] if $\bigl(\L_{\bar k},\O_{X_{\bar k}}(C)\bigr) = 0$ for a curve $C$ in $X_{\bar k}$ then $(C,C)_{X_{\bar k}} = (-2)$.
\end{enumerate}
\end{dfn}
If $(X,\l)$ is a polarized K3 surface over $k$, then one can find a finite separable extension $k'$ of $k$ such that $\l$ comes from a line bundle $\L_{k'}$ over $k'$. Indeed, this follows either from Remark \ref{EtPic} or from Proposition 4 in \cite[Ch. 8, \S 8.1]{NM} taking $T = \Spec(k^{\rm sp})$ and the fact that $\text{Br}(k^{\rm sp})$ is trivial.
\begin{dfn}\label{PolDefinition}
Let $S$ be scheme. A \emph{polarization} on a K3 space $\pi \colon X \lr S$ is a global section $\l\in \Pic_{X/S}(S)$ such that for every geometric point $\bar s$ of $S$ the section $\l_{\bar s} \in \Pic_{X_{\bar s}/\kappa(\bar s)}(\kappa(\bar s))$ is a polarization of $X_{\bar s}$. A \emph{quasi-polarization} on $X/S$ is a global section $\l\in \Pic_{X/S}(S)$ such that for every geometric point $\bar s$ of $S$ the section $\l_{\bar s} \in \Pic_{X_{\bar s}/\kappa(\bar s)}(\kappa(\bar s))$ is a quasi-polarization of $X_{\bar s}$.
\end{dfn}
\begin{dfn}
A polarization (respectively quasi-polarization) $\l$ on a K3 space $\pi \colon X \lr S$ is called \emph{primitive} if for every geometric point $\bar s$ of $S$ the polarization (respectively the quasi-polarization) $\l_{\bar s} \in \Pic_{X_{\bar s}/\kappa(\bar s)}(\kappa(\bar s))$ is primitive i.e., it is not a positive power of any element in  $\Pic_{X_{\bar s}/\kappa(\bar s)}(\kappa(\bar s))$.
\end{dfn}
\begin{lem}\label{PoltoLB}
Let $(\pi \colon X \lr S,\l)$ be a K3 space over $S$ with a polarization $\l$. Then one can find an \'etale covering $S' \lr S$ such that $\pi_{S'} \colon X_{S'} \lr S'$ is a K3 scheme and $\l_{S'}$ is the class of a relatively ample line bundle $\L_{S'}$ on $X_{S'}$.
\end{lem}
\begin{proof}
By definition one can find an \'etale covering $S_1 \lr S$ such that $\pi_1 \colon X_{S_1} \lr S_1$ is a K3 scheme. The pull-back $\l_{S_1}$ of $\l$ is a polarization on $X_{S_1}$. By Remark \ref{EtPic} the Picard functor $\Pic_{X_{S_1}/S_1}$ can be computed using the \'etale topology on $S_1$. Hence one can find an \'etale covering $S' \lr S$ such that $\l_{S'}$ is equal to the class of a line bundle $\L_{S'}$ on $X_{S'}$. By definition $\L_{S'}$ is pointwise ample hence using Lemma \ref{RelAmpl} we conclude that it is relatively ample. This finishes the proof.
\end{proof}
The self-intersection $(\L_{\bar s'},\L_{\bar s'})$ for a geometric point $\bar s'$ on $S'$ is constant on every connected component of $S'$. We say that $\l$ is a \emph{polarization of degree} $2d$ if $(\L_{\bar s'},\L_{\bar s'}) = 2d$ for every geometric point $\bar s'$ of $S'$.
\subsection{Automorphism Groups}
Let $S$ be a scheme and  $\pi \colon X \lr S $ be an algebraic space over $S$. Define the automorphism functor in the following way:
\begin{gather*}
 \Aut_S(X) \colon (\Sch/S)^0 \lr {\rm Groups} \\
\Aut_S(X)(T) = \Aut_T(X_T)
\end{gather*}
for every $S$-scheme $T$.
\begin{thm}\label{AutGrK3}
If $\pi \colon X \lr S$ is a polarized K3 space over $S$, then
$\Aut_S(X)$ is representable by a separated group scheme which is unramified and locally of finite type over $S$.
\end{thm}
\begin{proof}
Let $S' \lr S$ be an \'etale cover such that $\pi' \colon X' = X \times_S S' \lr S'$ is a projective K3 scheme over $S'$. The existence of such an \'etale covering $S'$ follows from Lemmas \ref{RelAmpl} and \ref{PoltoLB}. Let $S''$ be the product $S'\times_S S'$. Denote by $\pi_i$ the projection morphisms $\pi_i \colon X'\times_X X' \lr X' \lr X \lr S$ for $i = 1,2$. By definition $X'\times_X X'$ is representable by a quasi-compact subscheme of $X'\times_S X'$.

Using Proposition 1.4 in \cite[Ch. II]{Knu-AS} we can see that we have an exact sequence of groups
\begin{equation}\label{ReprAutGrSequence}
\xymatrix{
0 \ar[r] & \Aut_S(X)(T) \ar[r] & \Aut_{S'}(X')(T) \ar@<2pt>[r] \ar@<-2pt>[r] & \Aut_{S''}(X'\times_X X')(T).
}
\end{equation}
It follows from \cite[Exp. 221, \S4.c]{FGA} that the functors $\Aut_{S'}(X')$ and $\Aut_{S''}(X'\times_X X')$ are representable by group schemes locally of finite type over $S$. For simplicity we denote them by $\mathcal Y$ and $\mathcal W$ respectively. Then from the exact sequence \eqref{ReprAutGrSequence} we see that $\Aut_S(X)$ is representable by the fiber product
\begin{displaymath}
\xymatrix{
\Aut_S(X) \ar[d] \ar[rr] & & \mathcal W \ar[d]^{\Delta} \\
\mathcal Y \ar[rr]^{(pr_1^*,pr_2^*)} & & \mathcal W \times_S \mathcal W
}
\end{displaymath}
where $\Delta \colon \mathcal W \lr \mathcal W \times_S \mathcal W$ is the diagonal morphism.

The fact that the $\Aut_S(X)$ is separated follows directly from the valuative criterion for separatedness. 

To check that $\Aut_S(X)$ is unramified we may take $S$ to be the spectrum of an algebraically closed field $k$. A point in $\Aut_k(X)(k[\e]/(\e^2))$, which under the natural homomorphism maps to the the identity in $\Aut_k(X)(k)$, may be identified with a vector field on $X$. By Proposition \ref{HodgeNum} (1) a K3 surface has no non-trivial vector fields hence we conclude that $\Aut_k(X)$ is reduced.
\end{proof}
\begin{rem}
The proof of the theorem shows that $\Aut_S(X)$ is $0$-dimensional over $S$. Its fibers are constant group schemes.
\end{rem}
Let $\pi \colon X \lr S$ be a K3 space and let $\l$ be a polarization of $X$. Define the subfunctor $\Aut_S(X,\l)$ of $\Aut_S(X)$ in the following way
\begin{gather*}
 \Aut_S(X,\l) \colon (\Sch/S)^0 \lr {\rm Groups} \\
 \Aut_S(X,\l)(T) = \{ \a \in \Aut_S(X)(T)\ |\ \a^*\l= \l \in \Pic_{X/S}(T)\}
\end{gather*}
for every $S$-scheme $T$.
\begin{prp}\label{AutPolK3}
The functor $\Aut_S(X,\l)$ is a closed subfunctor of $\Aut_S(X)$. It is represented by a separated group scheme which is unramified and of finite type over $S$. Its relative dimension over $S$ is zero.
\end{prp}
\begin{proof}
The functor $\Aut_S(X,\l)$ is a closed subfunctor of $\Aut_S(X)$. It is representable by the subgroup scheme of $G = \Aut_S(X)$ (locally of finite type over $S$) given by the following (Cartesian) diagram:
$$
\xymatrix{
\Aut_S(X,\l) \ar[rr] \ar[d] & & S \ar[d]^{(\l,{\rm id})} \\
G = G\times_S S \ar[rr]^\psi & & \Pic_{X/S} \times_S S = \Pic_{X/S}.
}
$$
Here we have that $\l \colon S \lr \Pic_{X/S}$ is the section given by $\l$ and $\psi$ is the composition $\sigma \circ ({\rm id},\l)$ where 
$$
 \sigma \colon G \times \Pic_{X/S} \lr \Pic_{X/S}
$$
is the action of $G$ on $\Pic_{X/S}$.

Just as in the proof of the preceding theorem we may take $S$ to be the spectrum of an algebraically closed field $k$ in order to check that $\Aut_S(X,\l)$ is unramified. If $\a \in \Aut_k(X,\l)(k[\e]/\e^2)$ which is the identity in $\Aut_k(X,\l)(k)$, then by Theorem \ref{AutGrK3} above we see that $\a$ is the identity element of the group $\Aut_k(X)(k[\e]/\e^2)$. Since by definition we have an inclusion
$$
 \Aut_k(X, \l)(k[\e]/\e^2) \subset \Aut_k(X)(k[\e]/\e^2)
$$ 
we conclude that $\Aut_S(X,\l)$ is unramified over $S$.

Let $\bar s \colon \Spec(\Omega) \lr S$ be a geometric point. Then by \cite{Mat} (see also Corollary 2 in \cite{M-M}) the set $\Aut_S(X,\l)(\Omega)$ is finite. Hence $\Aut_S(X,\l)$ is of finite type over $S$.
\end{proof}
Note that in general, for a K3 surface $X$ over a field $k$, the group $\Aut_k(X)(k)$ might be infinite.
\begin{exa}
For any complex K3 surface $X$ with ${\rm rk}_\Z\Pic(X) = 20$ one has that $\Aut_\C(X)(\C)$ is infinite. For a proof see \cite[\S 5, Thm. 5]{Sh-I}.
\end{exa}
There are also examples of K3 surfaces $X$ having a finite group of automorphisms. An example of a complex K3 surface with ${\rm rk}_\Z\Pic(X) = 18$ and finite automorphism group is given in the remark on page 132 in \cite{Sh-I}.
\subsection{Automorphisms of Finite Order}In this section $k$ will be an algebraically closed field. If it is a field of characteristic $p$, then we will denote by $W$ the ring of Witt vectors with coefficients in $k$ and $K$ will be the field of fractions of $W$. 

Let $X$ be a  K3 surface over $k$. If $k= \C$, then it is a well-known theorem that $\Aut_\C(X)(\C)$ acts faithfully on $H^2_B(X,\Z)$. Here we prove a similar result for the automorphisms of finite order of $X$ acting trivially on $H^2_{\rm et}(X,\Z_l)$ where $l$ is a prime number different from ${\rm char}(k)$. The only restriction we impose is that ${\rm char}(k) \ne 2$. Later on in Section \ref{LevelStructures} we will introduce level structures on K3 surfaces and we will use this result to show that the corresponding moduli stacks are algebraic spaces.

\begin{lem}
Let $X$ be a K3 surface over $k$ and assume that {\rm char}$(k)=0$. Then $\Aut_k(X)(k)$ acts faithfully on $H^2_{\rm et}(X,\Z_l)$ for every prime $l$.
\end{lem}
\begin{proof}
Without loss of generality we may assume that the field $k$ can be embedded into $\C$. Fix an embedding $\sigma \colon k \hookrightarrow \C$. By the comparison theorem between Betti and \'etale cohomology we have an isomorphism $H^2_{\rm et}(X,\Z_l) \cong H^2_B(X \otimes_\sigma \C , \Z)\otimes_\Z \Z_l$. Let $\a \in \Aut_k(X)(k)$ be an automorphism acting trivially on $H^2_{\rm et}(X,\Z_l)$. Then $\a_\C$ acts trivially on $H^2_B(X \otimes_\sigma \C, \Z)\otimes \Z_l$. Since $H^2_B(X \otimes_\sigma \C, \Z)$ is a free $\Z$-module we conclude from \cite[Prop. 7.5]{L-P} that $\a = {\rm id}_X$.
\end{proof}
\begin{prp}\label{finaut}
Let $(X, \l)$ be a polarized K3 surface over $k$ and assume that $\cha(k) = p$ is different from $2$. Then the finite group $\Aut_k(X, \l)(k)$ acts faithfully on $H^2_{\rm et}(X,\Z_l)$ for any $l \ne p$.
\end{prp}
\begin{rem}
This result can be viewed as an analogue of Theorem 3 in \cite[Ch. IV]{Mum-AV} for (polarized) K3 surfaces.
\end{rem}
We will reduce the proof of Proposition \ref{finaut} to the preceding lemma. To do so we will use
crystalline cohomology and compare the action of an element in $\Aut_k(X,\l)(k)$ on $H^2_{\rm et}(X,\Q_l)$ and $H^2_{\rm cris}(X/K)$.

Let $X$ be a K3 surface over a field $k$. We denote by $H^n(X)$ and $H^n(X\times X)$ either $H^n_{\rm et}(X, \Q_l)$ and $H^n_{\rm et}(X \times X, \Q_l)$ for any $l$ prime to char$(k)$ or $H^n_{\rm cris}(X/K)$ and $H^n_{\rm cris}(X\times X/K)$. Note that we will be working with classes of certain algebraic cycles on $X$ and $X \times X$ so we should consider some Tate twists of these cohomology groups. But since $k$ is algebraically closed and the Galois action does not play any role in our consideration (we shall only consider some characteristic polynomials of automorphisms of $X$) we will omit these twists. 

For an isomorphism $\a \colon X \lr Y$ we will denote by $\a^*_l$ and $\a^*_{\rm cris}$ the isomorphisms induced on $H^2_{\rm et}(X,\Q_l)$ and $H^2_{\rm cris}(X/K)$ respectively. 
\begin{lem}\label{Kunncpts}
The K\"unneth components of the class $cl(u) \in H^4(X\times X)$ of any algebraic cycle on $X\times X$ are algebraic.
\end{lem}
\begin{proof} We have that $H^1_{\rm et}(X, \Q_l) = H^3_{\rm et}(X,\Q_l) = 0$ and $H^1_{\rm cris}(X/W) = H^3_{\rm cris}(X/W) = 0$. Then the K\"unneth
isomorphism reads
\begin{displaymath}
 H^4(X \times X) = \bigl(H^4(X)\otimes H^0(X)\bigr)\oplus \bigl(H^2(X)\otimes
 H^2(X)\bigr)\oplus\bigl(H^0(X)\otimes H^4(X)\bigr).
\end{displaymath}
Using this decomposition we write
\begin{displaymath}
  cl(u) = u_0 \oplus u_2 \oplus u_4.
\end{displaymath}
Every element of the one dimensional spaces $H^4(X)\otimes H^0(X)$
and $H^0(X)\otimes H^4(X)$ is algebraic. These are rational
multiple of the classes of $\{pt\}\times X$ and $X\times\{pt\}$.
Hence $u_0$ and $u_4$ are algebraic. It follows that $u_2$ is
expressed as a linear combination of algebraic classes, hence it is
algebraic.
\end{proof}
In particular, if $\Delta = \delta(X) \subset X\times X$ is the diagonal, then its K\"unneth components $cl(\Delta) = \pi_0 \oplus \pi_2 \oplus \pi_4 \in H^4(X\times X)$ are algebraic. Denote by $\langle \cdot , \cdot \rangle$ the intersection pairing on ${\rm CH}^2(X\times X)_\Q$.
\begin{cor}\label{ind}
Let $u \in {\rm CH}^2(X\times X)_\Q$ be a rational cycle and let $cl(u) \in
H^4(X\times X)$ be its algebraic class. Then its characteristic
polynomial $\det\bigl(1 - t\cdot cl(u)|H^2(X)\bigr)$ has rational coefficients
which are independent of $l$ and $p$ (i.e., of $H^2_{\rm et}(X,\Q_l)$ and $H^2_{\rm cris}(X/K)$). The coefficient in front of $t^i$ is given by
$$
 s_i = \langle u_i,\pi_2 \rangle
$$
for $i = 1,\dots,22$.
\end{cor}
\begin{proof}
The proof follows from the preceding lemma and by Theorem
3.1 in \cite{JT-Cycles}.
\end{proof}
\begin{thm}[Ogus]\label{Ogus}
If $p > 2$ then the natural morphism of groups
$$
 \Aut_k(X)(k) \lr \Aut\bigl(H^2_{\rm cris}(X/W)\bigr)
$$
is injective.
\end{thm}
\begin{proof} This is a result of A. Ogus and can be found in his paper on
Supersingular K3 crystals \cite[\S 2, Cor. 2.5]{O-SK3}.
\end{proof}
\noindent
{\it Proof of Proposition \ref{finaut}.} Take an element $\a \in \Aut_k(X,\l)(k)$. According to Proposition \ref{AutPolK3} it has finite order. Denote by $u =
\G_\sigma \subset X\times X$ the graph of $\a$. Then the automorphism of
$H^2_{\rm et}(X, \Q_l)$ induced by $cl(u) \in H^4_{\rm et}(X\times X, \Q_l)$ is the one induced by $\a$. By assumption it is the identity hence its characteristic polynomial is $(t-1)^{22}$. By Corollary \ref{ind} it is exactly the characteristic polynomial of the automorphism $\a^*_{\rm cris}$ of $H^2_{\rm cris}(X/K)$ induced by $\a$. Since $\a$ is an automorphism of finite order the induced map $\a^*_{\rm cirs}$ on the crystalline cohomology is semi-simple ($K$ has characteristic zero). Hence $\a^*_{\rm cris}$ acts trivially on $H^2_{\rm cris}(X/K)$ and by Theorem \ref{Ogus} it is the identity automorphism as $H^2_{\rm cris}(X/W)$ is torsion free.
\qed
\begin{rem}
Note that the only property of $\a$ which we used in the proof of Proposition \ref{finaut} is that it has finite order. This is really essential as in general the characteristic polynomial of $\a^*_l$ will not give enough information to conclude that the action of $\a_{\rm cris}$ on
$H^2_{\rm cris}(X/W)$ is trivial. The proof given above shows actually that any automorphism of finite order $\a$ of $X$ acting trivially on $H^2_{\rm et}(X, \Z_l)$ for some $l \ne p$ is the identity automorphism $\id_X$.
\end{rem}
\section{The Moduli Stack of Polarized K3 Surfaces}
We are ready to define moduli functors of (primitively) polarized K3 surfaces over $\Spec(\Z)$. We will follow the line of thoughts in \cite{D-M} in order to prove that these functors define Deligne-Mumford stacks. Shortly, this can be given in three steps.
\begin{enumerate}
\item [1.] Describe the deformations of primitively polarized K3 surfaces.
\item [2.] Construct a Hilbert scheme parameterizing K3 surfaces embedded in $\P^N$ for some appropriate $N \in \N$.
\item [3.] Construct a ``Hilbert morphism'' $\pi_{\rm Hib}$ from the Hilbert scheme to the moduli stack which is surjective and smooth. Use this morphism to conclude that the moduli stack is a Deligne-Mumford stack.  
\end{enumerate}
These steps are spelled out in detail in Sections \ref{DefTh}-\ref{ModStackSection}.
\subsection{Deformations of K3 Surfaces}\label{DefTh}
Let $k$ be an algebraically closed field. Denote by $W$
the ring of Witt vectors $W(k)$ in case $\cha(k) = p > 0$ and $W = k$ otherwise. Let $\underline A$ be the category of local artinian $W$-algebras $(A,{\mathfrak m}_A)$ together with an isomorphism $A/{\mathfrak m}_A \cong k$ compatible with the isomorphism $W/pW \cong k$. 

Let $X_0$ be a K3 surface over $k$. Consider the covariant functor
\begin{displaymath}
 \DefF_{\Sch}(X_0) \colon \underline A \lr {\rm Sets}
\end{displaymath}
given by
\begin{displaymath}
\begin{split}
  \DefF_{\Sch}(X_0)(A) =  \bigl\{ {\rm isom.\ classes\ of\ pairs}\ (X,\phi_0)\ |& {\rm where}\ X \lr \Spec(A) \\
   & {\rm is\ a\ K3\ scheme\ and}\ \phi_0\ {\rm is} \\
   & {\rm an\ isom.}\ \phi_0 \colon X\otimes_A k \cong X_0 \bigr\}. \\ 
\end{split}
\end{displaymath} 
\begin{prp}\label{DefFunct}
The functor $\DefF_{\Sch}(X_0)$ is pro-representable by a formal scheme $S$ over
$\Spf(W)$ which is formally smooth
of relative dimension $20$ i.e., it is (non-canonically) isomorphic to $\Spf(W[[t_1,\dots ,t_{20}]])$.
\end{prp}
\begin{proof} 
 This is Corollary 1.2 in \cite{Del-K3} in case $\cha(k) = p > 0$ and \cite[Cor. 5.7]{L-P} in case $\cha(k) =0$.
\end{proof}

Let $\L_0$ be a line bundle on $X_0$. For moduli problems one should study the deformations of 
the pair $(X_0,\L_0)$. Define
\begin{displaymath}
 \DefF_{\Sch}(X_0,\L_0) \colon \underline A \lr {\rm Sets}
\end{displaymath}
to be the functor sending an object $A$ of $\underline A$ to the isomorphism classes of triples $(X,\L, \phi_0)$ of flat deformations $X$ of $X_0$ over $A$, an invertible sheaf
$\L$ on $X$ and an isomorphism $\phi_0 \colon (X,\L)\otimes_A k \cong (X_0,\L_0)$. We have a morphism
\begin{equation}\label{DefSequence}
 \DefF_{\Sch}(X_0,\L_0) \lr \DefF_{\Sch}(X_0).
\end{equation}
\begin{thm}\label{poldef}
If the line bundle $\L_0$ is non-trivial, then the functor $\DefF_{\Sch}(X_0,\L_0)$ is pro-representable by a formally flat scheme of relative
dimension 19 over $W$ and the morphism \eqref{DefSequence} is a closed immersion, defined by a
single equation.
\end{thm}
\begin{proof}
See \cite[Prop. 1.5 and Thm. 1.6]{Del-K3}.
\end{proof}
Deligne proves that if $\L_0$ is an ample line bundle over $X_0$ then one can
find a discrete valuation ring $R$ which is a finite $W$ module and a lift
$(X \lr \Spec(R), \L) $ of $(X_0, \L_0)$ over $R$. In general one needs
ramified extensions of $W$ in order to find a lift of $(X_0, \L_0)$. The next
lemma shows that one can find a lift over $W$ if the self-intersection of $\L_0$ is prime to the characteristic of $k$. More precisely
one has:
\begin{lem}\label{smooth} 
Let $\L_0$ be an ample line bundle over $X_0$. If the polarization degree $(\L_0,\L_0)_{X_0} = 2d$ is prime to the characteristic of $k$, then $\DefF_{\Sch}(X_0,\L_0)$ is formally smooth.
\end{lem}
\begin{proof}
According to \cite[\S 2, Prop. 2.2 ]{O-SK3} (see also Lemma 2.2.6 in \cite{Del-CanCoord}) it is enough to see that $c_1(\L_0) \not \in F^2H^2_{DR}(X_0/k)$. Since we have that $(c_1(\L_0),c_1(\L_0)) = 2d \ne 0$ in $k$ it follows that $c_1(\L_0) \not \in F^2H^2_{DR}(X_0/k)$. For the proof in the case $k$ has characteristic zero we refer to \cite[\S 2, Thm. 1]{PSh-Sh}.
\end{proof}
\subsection{The Hilbert Scheme}\label{HilbertSchSection}
Recall that if $X$ is a K3 surface over a field $k$ with an ample line bundle $\L$, then the Hilbert polynomial of $\L$ is $h_\L(x)= dx^2+2$, where $(\L,\L) = 2d$. 

We fix two natural numbers $n$ and $d$ assuming that $n \geq 3$. Let $P_{d,n}(x)$ be the polynomial $n^2dx^2+2$ and let $N= P_{d,n}(1)-1$. Denote by ${\rm \bf Hilb}^{P_{d,n}}_N$ the \emph{Hilbert scheme} over $\Z$ representing the subvarieties of $\P^N$ with Hilbert polynomial $P_{d,n}(x)$. Let 
$$
\pi \colon \mathcal Z \lr {\rm \bf Hilb}^{P_{d,n}}_N
$$ 
be the universal family over the Hilbert scheme. For any morphism of schemes $f \colon S \lr {\rm \bf Hilb}^{P_{d,n}}_N$ we consider the following (Cartesian) diagram: 
\begin{equation}\label{HSchDiagram}
\xymatrix{
\mathcal X = S\times_{\bigl({\rm \bf Hilb}^{P_{d,n}}_N\bigr)} \mathcal Z \ar[rrr]^{f'} \ar[d]_{\pi'} & & & \mathcal Z \ar[d]_{\pi} \\
S \ar[rrr]^f & & & {\rm \bf Hilb}^{P_{d,n}}_N
}
\end{equation}
\begin{prp}\label{HilbSchK3}
There is a unique subscheme $H_{d,n}$ of ${\rm \bf Hilb}^{P_{d,n}}_N$ with the property:
\newline
A morphism of schemes $f \colon S \lr {\rm \bf Hilb}^{P_{d,n}}_N$ factors through $H_{d,n}$ if and only if the following conditions are satisfied.
\begin{enumerate}
\item [(i)] The pull-back $\X$ of the universal family over ${\rm \bf Hilb}^{P_{d,n}}_N$ is a K3 scheme over $S$ (see Diagram \eqref{HSchDiagram} above),
\item [(ii)] the line bundle ${f'}^* \O_{\P^N}(1)$ is isomorphic to $\L^n\otimes {\pi'}^*\M$ for some ample line bundle $\L$ on $\X$ and some line bundle $\M$ on $S$,
\item [(iii)] for every geometric point $\bar s \colon \Spec(\Omega) \lr S$ the natural homomorphism 
$$
 H^0(\P^N, \O_{\P^N}(1)) \otimes \Omega \lr H^0(\X_{\bar s}, \L^n_{\bar s})
$$
is an isomorphism.
\end{enumerate} 
There exists an open subscheme $H^{pr}_{d,n}$ of $H_{d,n}$ such that: A morphism of schemes $f \colon S \lr {\rm \bf Hilb}^{P_{d,n}}_N$ factors through $H_{d,n}^{pr}$ if and only if conditions {\it (i), (ii)} and {\it (iii)} are satisfied and in addition for every geometric point $\bar s$ of $S$ the line bundle $\L_{\bar s}$ from {\it (ii)} is primitive.
\end{prp}
\begin{proof}
The proof of the proposition is standard and can be found in the case of curves in Mumford's book \cite[Ch. 5, \S 2, Prop. 5.1]{Mum-GIT}. We shall sketch only the additional arguments needed in our situation. 

There is a maximal open subscheme $U_1$ of ${\rm \bf Hilb}^{P_{d,n}}_N$ such that every fiber of the pull-back $\X_1$ of the universal family $\mathcal Z$ over $U_1$ is a non-singular variety. Let $U_2$ be the open subscheme of $U_1$ consisting of the points $s$ for which $H^1(\X_{1,s}, \O_{\X_{1,s}}) = 0$ (see \cite[Ch. III, \S 12, Thm. 12.8]{HAG}). Denote by $\X_2$ the pull-back of the universal family over $U_2$.

Let $\Pic_{\X_2/U_2}$ be the relative Picard scheme of $\X_2$ over $U_2$. The two line bundles $\Omega^2_{\X_2/U_2}$ and $\O_{\X_2}$ define two morphisms: $\omega, \l \colon U_2 \lr \Pic_{\X_2/U_2}$. Define $U_2$ to be the fiber product:
$$
\xymatrix{
U_3 \ar[d] \ar[rrr] & & & U_2 \ar[d]^{(\l,\omega)} \\
\Pic_{\X_2/U_2} \ar[rrr]^{\Delta} & & & \Pic_{\X_2/U_2} \times_{U_2} \Pic_{\X_2/U_2}
}
$$ 
where $\Delta$ is the diagonal morphism. Since $\Pic_{\X_2/U_2}$ is separated $U_3$ is a closed subscheme of $U_2$. The pull-back $\X_3 \lr U_3$ of the universal family over ${\rm \bf Hilb}^{P_{d,n}}_N$ is a K3 scheme.

Let $[n] \colon \Pic_{\X_3/U_3} \lr \Pic_{\X_3/U_3}$ be the multiplication by-$n$-morphism. The pull-back of $\O_{\P^N}(1)$ over $U_3$ defines a morphism $\l \colon U_3\lr \Pic_{\X_3/U_3}$. Define $U_4$ to be the fiber product
$$
\xymatrix{
U_4 \ar[d] \ar[rr] & & U_3 \ar[d]^\l \\
\Pic_{\X_3/U_3} \ar[rr]^{[n]} & & \Pic_{\X_3/U_3}.
}
$$
By Lemma \ref{PicDiv} the morphism $[n]$ is a closed immersion hence $U_4$ is a closed subscheme of $U_3$. Clearly, $U_4$ is the subscheme of ${\rm \bf Hilb}^{P_{d,n}}_N$ for which properties (i) and (ii) hold. One takes $H_{d,n}$ to be the (closed) subscheme of $U_4$ obtained as in the end of the proof of Proposition 5.1 in \cite[Ch. 5, \S 2]{Mum-GIT} (where instead of $\Omega^1_{\Gamma/U_2}$ one works with the pull-back $\L'$ of the bundle $\O_{\P^N}(1)$). It satisfies all conditions of the proposition. 

To show the existence of $H^{pr}_{n,d}$ one has to take the open subscheme $U^0_4$ of $U_4$ above corresponding to the points in $U_4$ over which the class of the pull-back of $\O_{\P^N}(1)$ in $\Pic_{\X_4/U_4}$ is only divisible by $n$. The existence of such a subscheme can be seen, in a way similar to the proof of Lemma \ref{PicDiv} {\it (ii)}, using the fact that the homomorphisms $[n] \colon \Pic_{\X_4/U_4} \lr \Pic_{\X_4/U_4}$ are closed immersions. 
\end{proof}
We will use the schemes $H_{d,n}$ and $H^{pr}_{d,n}$ to construct moduli stacks of polarized K3 surfaces over $\Z$.
\subsection{The Moduli Stack}\label{ModStackSection}
One way to construct the coarse moduli space of complex K3 surfaces with a primitive
polarization of degree $2d$ is to use period maps. This approach is taken up in \cite[Expos\'e XIII, \S 3]{Ast-K3}. Here we will use rather
different techniques to deal with this problem in positive and more generally in mixed
characteristic.
 
\begin{dfn}\label{modsp}
Let $d$ be a natural number. Consider the category $\Fk_{2d}$ defined in the following way:
\begin{enumerate}
\item[{\it Ob:}] The objects of $\Fk_{2d}$ are pairs $(\pi \colon X \lr S, \l)$ consisting of a K3 space $\pi \colon X \lr S$ with a primitive polarization $\l$ of degree $2d$ over $S \in \Sch$.
\item [{\it Mor:}] For two objects $\X_1 = (\pi_1 \colon X_1 \lr S_1, \l_1)$ and
$\X_2 = (\pi_2 \colon X_2 \lr S_2, \l_2)$ we define the morphisms to be
 \begin{displaymath}
 \begin{split}
 \Hm(\X_1, \X_2)  = \bigl\{{\rm pairs}\ (f_S,f)\ |& f_S \colon S_1 \lr S_2 \ {\rm is\ a\
 morph.\ of} \\
  & {\rm schemes\ and}\ f \colon X_1 \lr X_2\times_{S_2,f_S} S_1 \\
  & {\rm is\ an\ isom.\ over}\ S_1\ {\rm with}\ f^*\l_2 = \l_1\bigr\}. \\
\end{split}
\end{displaymath}
\end{enumerate}
\end{dfn}
The functor $p_{\Fk_{2d}} \colon \Fk_{2d} \lr \Sch$ sending
a pair $(\pi \colon X \lr S, \l)$ to $S$ makes $\Fk_{2d}$ into a category over
$\Sch$. We will denote by $\Fk_{2d,S}$ the full subcategory
of $\Fk_{2d}$ consisting of the objects over $S$.
\begin{dfn}
For a natural number $d$ we define the category $\mathcal M_{2d}$ of K3 spaces with a polarization of degree $2d$ in the same way as in Definition \ref{modsp} but taking as objects pairs of polarized K3 spaces $(\pi \colon X \lr S, \l)$ over a scheme $S$.
\end{dfn}
We have that $\Fk_{2d}$ is a full subcategory of $\mathcal M_{2d}$. Those two categories are the same if and only if $d$ is square-free. 
\begin{thm}\label{F-Prp}
 The categories $\Fk_{2d}$ and $\mathcal M_{2d}$ are separated
Deligne-Mumford stacks of finite type over $\Z$. The inclusion $\Fk_{2d} \hookrightarrow \mathcal M_{2d}$ is an open immersion.
\end{thm}
\begin{dfn}
We will call $\Fk_{2d}$ \emph{the moduli stack of primitively polarized K3 surfaces of degree} $2d$ and $\mathcal M_{2d}$ \emph{the moduli stack of polarized K3 surfaces of degree} $2d$.
\end{dfn}
\begin{rem}
Let us explain first why we want to consider moduli of \emph{primitively} polarized K3 surfaces. For various reasons we will have to work with algebraic spaces rather than with algebraic stacks. Just like in the case of abelian varieties one can introduce level structures on K3 surfaces and hope that the corresponding moduli problems are representable by algebraic spaces. We will define level structures on a polarized K3 surface $(X,\l)$ using its primitive cohomology groups $P^2_{\rm et}(X_{\bar k},\Z_l(1))$ for certain primes $l$ (see Section \ref{LevelStructures}). To be able to do that we will need that $P^2_{\rm et}(X_{\bar k},\Z_l(1))$ belongs to a single isometry class of quadratic lattices, which is the case, if $\l$ is primitive.
\end{rem}
We will prove the theorem in a sequence of steps.
\begin{lem}
The categories $\Fk_{2d}$ and $\mathcal M_{2d}$ are groupoids.
\end{lem}
\begin{proof}
We have to check two axioms. See for instance \cite[Ch. 2, Def. 2.1]{L-MB} or 
p. 96 of \cite{D-M}. One sees immediately that the usual notions of pull-backs satisfy these two axioms.
\end{proof}
\begin{lem}\label{LStack}
The groupoids $\Fk_{2d}$ and $\mathcal M_{2d}$ are stacks for the \'etale topology.
\end{lem}
\begin{proof}
The proofs for $\M_{2d}$ and $\Fk_{2d}$ are exactly the same so we will prove the lemma for $\Fk_{2d}$. We have to check two properties. Namely, first we will show that for any scheme $S \in \Sch$ and any two
objects $\X$ and $\mathcal Y$ over $S$ the functor
\begin{displaymath}
 \Is_S(\X,\mathcal Y): (\Sch/S) \lr {\rm Sets}
\end{displaymath}
defined by
$$
 (\pi \colon S' \lr S) \mapsto \Hm(\pi^*\X,\pi^*\mathcal Y)
$$
is a sheaf for the \'etale topology on $S$. Then we prove that descent data are
effective (cf. \cite[Ch. 2, Def. 3.1]{L-MB} or Definition 4.1 in \cite{D-M}). 
\newline
\newline
{\it The functor} $\Is_S(\X,\mathcal Y)$ {\it is an \'etale sheaf:} Take two objects $\X = (X \lr S, \l_{X})$ and $\mathcal Y = (Y \lr S, \l_Y)$ over $S$. Let $S'$ be an $S$-scheme.

Let $\{S'_i\}_{i \in I}$ be an \'etale covering of $S'$ and $f_j \in \Is_S(\X,\mathcal Y)(S')$ for $j = 1,2$ are two elements such that
$f_1|_{S'_i} = f_2|_{S'_i}$. Then clearly $f_1 = f_2$ as isomorphisms of the
pair $(\X_{S'}, \mathcal Y_{S'})$.

Let $\{S'_i\}_{i \in I}$ be an \'etale covering of $S'$. Suppose given elements $f_i \in \Is_S(\X,\mathcal Y)(S'_i)$ such that $f_i|_{S'_{ij}} =
f_j|_{S'_{ij}}$ where $S'_{ij} = S'_i\times_{S'} S'_j$. We have to show that 
those come from a global ``isomorphism''. Note that without loss of generality we may assume that $X_i \lr S'_i$ are K3 schemes. Combining \cite[Ch. II, Prop. 1.4]{Knu-AS} and effectiveness of descent for morphisms of schemes (see
\cite[Ch. 6, \S 1, Thm. 6(a)]{NM}) we conclude that $f'$ descends to a morphism $f
\colon X_{S'} \lr Y_{S'}$ such that $f_{S_i'} = f_i$.
Since $\Pic_{X/S}$ and $\Pic_{Y/S}$ are algebraic spaces (in particular sheaves for the \'etale topology on $S$) and
$f^*\l_{Y_{S'}}|_{S_i'} = \l_{X_{S'}}|_{S_i'}$ we see that $f^*\l_{Y_{S'}} = \l_{X_{S'}}$. Hence we have that $f \in \Is_S(\X,\mathcal Y)(S')$ and $f|_{S_i'} = f_i$. This
shows that $\Is_S(\X,\mathcal Y)$ is an \'etale sheaf.
\newline
\newline
{\it Effectiveness of descent:} Suppose given an \'etale cover $S'$ of $S$
and an object $\X' = (\pi' \colon X' \lr S', \l')$ with descent datum over $S$. Without loss of generality we may assume that the algebraic space $X'$ is actually a scheme (by refining the \'etale covering $S'$ if needed). We have to show that $(\pi' \colon X' \lr S',\l')$ descends to a polarized K3 space $(\pi \colon X \lr S,\l)$ over $S$. 

Denote by $S''$ the product $S' \times_S S'$ and let $pr_i$ for $i=1,2$ be the two projection maps. The descent datum on $X' \lr S'$ over $S$ identifies the two schemes $pr_1^*X'$ and $pr_2^*X'$. Denote this scheme by $R$. Then we have two \'etale morphisms
\begin{displaymath}
\xymatrix{ 
 R \ar@<2pt>[r] \ar@<-2pt>[r] & X'
}
\end{displaymath} 
which make $R \subset X'\times_S X'$ into an \'etale equivalence relation. Following the constructions of \cite[Ch. I, \S 5, 5.4]{Knu-AS} we obtain an algebraic space $X$ over $S$ such that $X\times_S S'$ is isomorphic to $X'$. Hence $\pi \colon X \lr S$ is a K3 space. 

Since $\Pic_{X/S}$ is an \'etale sheaf the local section $\l'$ over $S'$ together with descent datum over $S$ give rise to a global section $\l \in \Pic_{X/S}(S)$ such that $\l_{S'} = \l'$. Clearly, $\l$ is a polarization of $X \lr S$.
\end{proof}
Next we deal with the representability of the isomorphism functors of
polarized K3 surfaces. For two algebraic spaces $X$ and $Y$ over a base scheme $S$ define the contravariant isomorphism functor
$$
 \Is_S(X,Y) \colon (\Sch/S) \lr {\rm Sets}
$$ 
by
$$
 \Is_S(X,Y)(T) = \{f \colon X_T \lr Y_T\ |\ f\ {\rm is\ an\ isomorph.\ of\ alg.\ spaces\ over}\ T\}
$$ 
for any $S$-scheme $T$.
\begin{lem}\label{Isom}
For any $S \in \Sch$ and two objects $\X$ and $\mathcal Y$ of $\Fk_{2d}$ (respectively $\mathcal M_{2d}$) over $S$, the
functor $\Is_S(\X,\mathcal Y)$ is representable by a separated scheme which is unramified and of finite type over $S$.
\end{lem}
\begin{proof}
Let $\X$ and $\mathcal Y$ be the objects $(X \lr S,\l_X)$ and $(Y \lr S, \l_Y)$ respectively. 
\newline
\newline
{\bf Step 1:}
We can find an \'etale cover $S'$ of $S$ such that $X' = X\times_S S'$ and $Y' = Y\times_S S'$ are projective K3 schemes over $S'$. Denote by $S''$ the product $S'\times_S S'$. By \cite[Exp. 221, \S4.c]{FGA}) the functors $\Is_{S'}(X',Y')$ and $\Is_{S''}(X'\times_X X', Y'\times_Y Y')$ are representable by schemes $\mathcal U$ and $\mathcal V$, locally of finite type over $S$. By Proposition 1.4 in \cite[Ch. II]{Knu-AS} one has an exact sequence of sets
\begin{displaymath}
\xymatrix{
0 \ar[r] & \Is_S(X,Y)(T) \ar[r] & \Is_{S'}(X',Y')(T) \ar@<2pt>[r] \ar@<-2pt>[r] & \Is_{S''}(X'\times_X X', Y'\times_Y Y')(T).
}
\end{displaymath}
Then we see that $\Is_S(X,Y)$ is representable by the scheme defined by the following Cartesian diagram
\begin{displaymath}
\xymatrix{
\Is_S(X,Y) \ar[d] \ar[r] & \mathcal V \ar[d]^{\Delta} \\
\mathcal U \ar[r]^{(pr_1^*,pr_2^*)} & \mathcal V \times_S \mathcal V
}
\end{displaymath}
where $\Delta \colon \mathcal V \lr \mathcal V \times_S \mathcal V$ is the diagonal morphism. 
\newline
\newline
{\bf Step 2:} 
By Step 1 the functor $\Is_S(X,Y)$ is represented by a scheme locally of finite type over $S$. Then the functor $\Is_S(\X,\mathcal Y)$ is represented by the scheme defined by the following Cartesian diagram:
$$
\xymatrix{
\Is_S(\X,\mathcal Y) \ar[rrr] \ar[d] & & & S \ar[d]^{\l_X} \\
\Is_S(X,Y) \times_S S \ar[rr]^{({\rm id},\l_Y)} & &
\Is_S(X,Y) \times_S \Pic_{Y/S} \ar[r] & \Pic_{X/S}
}
$$
where the bottom-right arrow is just the pull back morphism.
\newline
\newline
{\bf Step 3:}
We are left to show that $\Is_S(\X,\mathcal Y)$ is unramified over $S$. As in the proof of Theorem \ref{AutGrK3} it is enough to check the properties of $\Is_S(\X,\mathcal Y)$ when $S$ is a spectrum of an algebraically closed field. In this case $\Is_S(\X,\mathcal Y)$ is either empty or it is isomorphic to $\Aut_k(X,\l)$. As the latter is separated, reduced and of finite type over $k$ we conclude that the same holds for $\Is_S(\X,\mathcal Y)$.
\end{proof}
\noindent
{\it Proof of Theorem \ref{F-Prp}:} We will give the proof for $\Fk_{2d}$ in several steps. For the proof that $\mathcal M_{2d}$ is a Deligne-Mumford stack one should only replace $H^{pr}_{d,3}$ by $H_{d,3}$ below.
\newline
\newline
{\bf Step 1:} We saw in Proposition \ref{HilbSchK3} that there exists a Hilbert scheme $H_{3,d}^{pr}$, of finite type over $\Z$, classifying K3 surfaces with a polarization of degree $2d$ which are embedded in a projective space via the third power of the polarization. One has then the universal family $f \colon \X \lr H^{pr}_{3,d}$ and we know that $\O_\X(1) \cong \L^3 \otimes f^*\M$ for some ample line bundle $\L$ on $\X$ of degree $2d$ and an invertible sheaf $\M$ on $H^{pr}_{3,d}$. Although the line bundle $\L$ with this property is not unique, its class $\l_\X = cl(\L) \in \Pic_{\X/H^{pr}_{3,d}}$ is uniquely determined as $\l_\X^3 = cl(\O_\X(1))$. Define the morphism of stacks
\begin{displaymath}
 \pi_{\rm Hilb} \colon H_{3,d}^{pr} \lr \Fk_{2d}.
\end{displaymath}
sending $H^{pr}_{3,d}$ to the pair $(f\colon \X \lr H^{pr}_{3,d},\l_\X)$. By construction the self-intersection index $(\l_{\X,h},\l_{\X, h})$ is $2d$ for any $h \in H^{pr}_{3,d}$ and $\l_\X$ is primitive so this morphism is correctly defined. 
\newline
\newline
{\bf Step 2:} {\it The morphism $\pi_{\rm Hilb}$ is surjective}. This follows form the definition (cf. \cite[Def. 3.6]{L-MB}) and Lemma \ref{PoltoLB}. Indeed, for any $(\pi \colon X \lr S, \l) \in \Fk_{2d}(S)$ one can find an \'etale cover $S' \lr S$ such that $\pi_{S'} \colon X_{S'} \lr S'$ is a K3 scheme and $\l_{S'}$ is equal to the class of a relatively ample line bundle $\L'$ on $X_{S'}$. By Lemma \ref{RelAmpl} the line bundle ${\L'}^3$ defines a closed immersion $X_{S'} \hookrightarrow \P({\pi_{S'}}_*{\L'}^3)$. Refining $S'$ further if needed we may assume that $\P({\pi_{S'}}_*{\L'}^3)$ is isomorphic with $\P^{9d+1}_{S'}$. Then the inclusion $X_{S'} \hookrightarrow \P({\pi_{S'}}_*{\L'}^3)$ satisfies the conditions of Proposition \ref{HilbSchK3} by construction. Hence it corresponds to a morphism $f_X \colon S' \lr H^{pr}_{3,d}$ and we have that
\begin{displaymath}
 \pi_{\rm Hilb}(f_X \colon S' \lr H^{pr}_{3,d}) = (\pi_{S'} \colon X_{S'} \lr S', \l_{S'}).
\end{displaymath}
\newline
\newline
{\bf Step 3:} {\it The morphism $\pi_{\rm Hilb}$ is representable and smooth}. Let $S$ be a scheme and suppose given a morphism $S \lr \Fk_{2d}$ corresponding to a primitively polarized K3 space $(\pi \colon X \lr S, \l)$. We have to show that the product $S\times_{\Fk_{2d}} H^{pr}_{3d}$ is representable by an algebraic space which is smooth over $S$ (via $pr_1$). By the surjectivity of $\pi_{\rm Hilb}$ one can find an \'etale cover $S'$ of $S$ and a projective embedding $X_{S'} \hookrightarrow \P^{9d+1}_{S'}$, defined by a very ample line bundle $\L^3$. It gives rise to a morphism $S' \lr H^{pr}_{3,d}$ with 
$$
 \pi_{\rm Hilb}(S' \lr H^{pr}_{3,d}) = (X_{S'} \lr S',\l_{S'}) \in \Fk_{2d}(S').
$$
We claim that the product $S'\times_{\Fk_{2d}} H^{pr}_{3,d}$ is representable by a scheme isomorphic to $\PrG(9d+2)_{S'}$. For any $S'$-scheme $U$ we have that 
\begin{displaymath}
\begin{split}
S'\times_{\Fk_{2d}} H^{pr}_{3,d}(U)  = &  \biggr\{\ \bigr((U \lr S'),(U \lr H^{pr}_{3,d}),g\bigl)\ \ \big|\\ 
& \ g\in \Hm\bigr((X_U \lr U, \l_U),\pi_{\rm Hilb}(U \lr H^{pr}_{3,d})\bigl)\ {\rm in}\ \Fk_{2d}\ \biggl\} \\
\end{split}
\end{displaymath}
where $\pi_{\rm Hilb}(U \lr H^{pr}_{3,d}) = (\X_U \lr U, \l_\X|_U)$. Any such morphism $g$ gives rise to an isomorphism $\L^3 \cong \O_{\X_U}(1)\otimes f_U^*\M$ for some invertible sheaf $\M$ on $U$ and hence an isomorphism 
$$
 \P(\pi_{U*}\L^3) \cong \P\bigl(f_{U*}\O_{\X_U}(1)\otimes \M\bigl).
$$
But by condition $(iii)$ of Proposition \ref{HilbSchK3} we have an isomorphism 
$$
 \P(f_{U*}\O_{\X_U}(1)\otimes \M) \cong \P\bigr(pr_{U*}\O_{\P^{9d+1}_U}(1)\bigr) = \P^{9d+1}_U
$$ 
and hence we obtain an isomorphism $\P(\pi_{U*}\L^3) \cong \P^{9d+1}_U$. This correspondence gives a bijection 
\begin{displaymath}
S'\times_{\Fk_{2d}} H^{pr}_{3,d}(U) \leftrightarrow \bigl\{\ {\rm isomorphisms}\ \P(\pi_{U*}\L^3) \cong \P^{9d+1}_U\ \bigr\}
\end{displaymath} 
and the right hand set can be identified with $\PrG(9d+1)_{S'}(U)$. For this we refer to the arguments given on pp. 101-103 in \cite{Mum-GIT}. Hence $S'\times_{\Fk_{2d}} H^{pr}_{3,d}$ is representable by the scheme $\PrG(9d+1)_{S'}$ which is smooth over $S'$.

We will show next that $S\times_{\Fk_{2d}} H^{pr}_{3,d}$ is a smooth algebraic space over $S$. We have a surjective map of \'etale sheaves 
$$
 S'\times_{\Fk_{2d}} H^{pr}_{3,d} \lr S\times_{\Fk_{2d}} H^{pr}_{3,d}.
$$
The product
\begin{displaymath}
R := \bigl(S'\times_{\Fk_{2d}} H^{pr}_{3,d}\bigr) \times_{\bigl(S\times_{\Fk_{2d}} H^{pr}_{3,d}\bigr)} \bigl(S'\times_{\Fk_{2d}} H^{pr}_{3,d}\bigr)
\end{displaymath}
can be identified with the smooth $S$-scheme $(S'\times_S S')\times_{\Fk_{2d}} H^{pr}_{3,d}$. The natural morphism 
$$
 R \lr (S'\times_{\Fk_{2d}} H^{pr}_{3,d}) \times (S'\times_{\Fk_{2d}} H^{pr}_{3,d})
$$
is quasi-compact and the two projection maps
\begin{displaymath}
\xymatrix{ 
 R \ar@<2pt>[r] \ar@<-2pt>[r] & S\times_{\Fk_{2d}} H^{pr}_{3,d}
}
\end{displaymath} 
are \'etale as they correspond to the two \'etale projection morphisms $\xymatrix{S'\times_S S' \ar@<2pt>[r] \ar@<-2pt>[r] & S}$. Hence $S\times_{\Fk_{2d}} H^{pr}_{3,d}$ is an algebraic space, which is moreover smooth over $S$ as it possesses a smooth atlas $S'\times_{\Fk_{2d}} H^{pr}_{3,d}$ (over $S$).
\newline
\newline
{\bf Step 4:} Using Remark 4.1.2 (i) in \cite[Ch. 4]{L-MB} (or Prop. 4.4 in \cite{D-M}) and Lemma \ref{Isom} we see that the diagonal morphism $\Delta \colon \Fk_{2d} \lr \Fk_{2d} \times_\Z \Fk_{2d}$ is representable, separated and quasi-compact. Then we can apply Theorem 4.21 of \cite{D-M} to the morphism $\pi_{\rm Hilb} \colon H_{3,d}^{pr} \lr \Fk_{2d}$ and conclude that $\Fk_{2d}$ is a Deligne-Mumford stack of finite type over $\Z$.
\newline
\newline
{\bf Step 5:} We will show that the algebraic stack $\Fk_{2d}$ is separated. 
As $\Fk_{2d}$ is of finite type over $\Z$ one can use the valuative criterion for separateness from \cite[Thm. 4.18]{D-M} (cf. \cite[Prop. 7.8 and Thm. 7.10]{L-MB}). It reduces to showing that if $(\pi_i \colon X_i \lr S,
\l_i)$, for $i= 1,2$, are two primitively polarized K3 spaces over the spectrum $S$ of a discrete 
valuation ring $R$ with field of fractions $K$, then every isomorphism 
$f \colon (X_1 \otimes K, \l_1 \otimes K) \lr  (X_2 \otimes K, \l_2 \otimes K)$ extends
to a $S$-isomorphism between $(X_1, \l_1)$ and  $(X_2, \l_2)$. Note that after taking a finite \'etale covering of $S$ we may assume that:
\begin{enumerate}
\item [(a)] $X_i$ are schemes,
\item [(b)] $\l_i = c_1(\L_i)$ for some ample line bundle $\L_i$,
\item [(c)] $f$ gives an isomorphism of pairs $f \colon  (X_1 \otimes K, \L_1 \otimes K) \lr (X_2 \otimes K, \L_2, \otimes K)$. 
\end{enumerate}
Then using \cite[Thm. 2]{M-M} (as a K3 surface is non-ruled) we see that $f$ extends uniquely to an isomorphism between $(X_1, \L_1)$ and $(X_2, \L_2)$.
\newline
\newline
{\bf Step 6:} We are left to show that the natural inclusion $\Fk_{2d} \hookrightarrow \mathcal M_{2d}$ is an open immersion. Take a noetherian scheme $S$ and suppose given a morphism $S \lr \mathcal M_{2d}$ corresponding to a polarized K3 space $(\pi \colon X \lr S, \l)$. Let $f\colon S' \lr S$ be an \'etale covering such that $\pi_{S'} \colon X_{S'} \lr S'$ is strongly projective (cf. Step 2 in the proof of Theorem \ref{F-Prp}). According to Lemma \ref{PicDiv} the set of points 
$$
 {S'}^o = \{s \in S'\ |\ {\rm such\ that}\ \l_{S',s}\ {\rm is\ primitive}\}
$$ 
is an open subscheme of $S'$. The morphism $f$ is \'etale and hence $f({S'}^o) \subset S$ is also an open subscheme which represents $S\times_{\mathcal M_{2d}} \Fk_{2d}$.
\qed
\begin{rem}
Another possible proof of Theorem \ref{F-Prp} is to use Artin's criterion (\cite[Cor. 10.11]{L-MB}). This approach is taken up in \cite[Thm. 6.2]{Olsson-K3} where M. Olsson constructs a compact stack of ``polarized log K3 spaces'' over $\Q$. 
\end{rem}
An immediate consequence of Theorem \ref{F-Prp} is the existence of a coarse moduli space of polarized K3 surfaces.
More precisely Corollary 1.3 in \cite{K-M} says
\begin{cor}\label{CoarseSpaces}
The moduli stacks $\Fk_{2d}$ and $\mathcal M_{2d}$ have coarse moduli spaces which are separated algebraic spaces.
\end{cor}
Note that this argumet shows that $\Fk_{2d}$ and $\mathcal M_{2d}$ are global quotient stacks.

Before going on we will shortly outline how one can obtain stronger results on coarse moduli schemes of polarized K3 surfaces in characteristic zero.
\newline
\newline
{\bf Approach via periods of K3 surfaces.}\label{ModSpPeriodAppr} As we mentioned in the beginning of this section one can use analytic methods to construct a coarse moduli scheme of primitively polarized K3 surfaces. Consider the complex space
\begin{displaymath}
\Omega^\pm = \{ \omega \in \P(L_{2d}\otimes \C) |\ \psi_{2d}(\omega,\omega) = 0\ \text{and}\ \psi_{2d}(\omega,\bar \omega) > 0 \}
\end{displaymath} 
which consists of two connected components. It can be identified with the space 
$$
 \SO(2,19)(\R)/\bigl(\SO(2)(\R) \times \SO(19)(\R)\bigr).
$$
Let $\Omega^+$ denote one of its connected components, say corresponding to 
$$
\SO(2,19)(\R)^+/\bigl(\SO(2)(\R) \times \SO(19)(\R)\bigr),
$$
where $\SO(2,19)(\R)^+$ is the connected component of $\SO(2,19)(\R)$ containing the identity. It is a bounded symmetric domain of type IV and of dimension 19. Let $\Gamma$ be the group $\{g \in {\rm O}(V_0)(\Z)\ |\ g(e_1-df_1) = e_1-df_1\}$ and denote by $\Gamma^+$ the subgroup of $\Gamma$ of index 2 which consists of isometries preserving the connected components of $\Omega^\pm$. Then $\Gamma^+$ acts on $\Omega^+$ properly discontinuously and the space $\Omega^+/\Gamma^+$ is a coarse moduli scheme for primitively quasi-polarized complex K3 surfaces of degree $2d$. There is an open part $\Omega^0$ of $\Omega^+$ such that $\Omega^0/\Gamma^+$ is a coarse moduli scheme for primitively polarized complex K3 surfaces of degree $2d$. For details and proofs we refer to \cite[Exp. XIII]{Ast-K3}. The existence of a coarse moduli scheme is Proposition 8 in {\it loc. cit.}.
\newline
\newline
{\bf Approach via geometric invariant theory.} Let $k$ be an algebraically closed field of characteristic zero. Then using the techniques of \cite[Ch. 8]{Vie-M}, and more precisely \S 8.2 (see Theorem 8.23), one can prove that the moduli functor $\Fk_{2d}\otimes k$ (respectively $\mathcal M_{2d} \otimes k$) has a quasi-projective coarse moduli scheme over $k$. Indeed, one has that Assumptions 8.22 in \cite[\S 8.2]{Vie-M} are satisfied:
\begin{enumerate}
\item [(i)] The functor is locally closed. This follows from the proof of Proposition \ref{HilbSchK3}. 
\item [(ii)] The separateness property is shown in Step 2 of the proof of Theorem \ref{F-Prp}. 
\item [(iii)] The functor is bounded by Theorem \ref{albprop}. See also Remark 8.24 in {\it loc. cit.} and note that the condition `$\omega^2$ is trivial' is a locally closed condition. 
\end{enumerate}
One actually shows that the scheme in question is $H_{3,d}^{pr}\otimes k/ \PrG (N)_k$ (respectively $H_{3,d}\otimes k/ \PrG (N)_k$) for a suitable $N \in \N$.

Combining the approach to coarse moduli schemes via geometric invariant theory and Corollary \ref{CoarseSpaces} we conclude that $\Fk_{2d, \Q}$ (respectively $\mathcal M_{2d,\Q}$) has a quasi-projective coarse moduli scheme.
\begin{prp}\label{smoothstac}
The moduli stacks $\Fk_{2d}$ and $\mathcal M_{2d}$ are smooth of relative dimension 19 over $\Z[\frac{1}{2d}]$.
\end{prp}
\begin{proof}
According to \cite[Prop. 4.15]{L-MB} we have to show that for any strictly henselian local ring $R$ and surjection $\Spec(R) \lr \Spec(R')$ defined by a nilpotent sheaf of ideals one has that the natural map
$$
 {\rm Hom}(\Spec(R'), \Fk_{2d,\Z[1/2d]}) \lr {\rm Hom}(\Spec(R),\Fk_{2d,\Z[1/2d]}) 
$$ 
is surjective. Since $R$ is strictly henselian every K3 space over $\Spec(R)$ is a K3 scheme and the same holds for spaces over $\Spec(R')$ (see \cite[EGA IV, 18.1.2]{EGA}). Hence by Lemma \ref{smooth} we conclude that $\Fk_{2d,\Z[1/2d]}$ is smooth over $\Z[1/2d]$.

The same argument applies also to the dimension claim. Since every K3 space over $\Spec(k[\e]/\e^2)$ is a K3 scheme we conclude from Theorem \ref{poldef} that the dimension of $\Fk_{2d,\Z[1/2d]}$ at every point is 19.

This proof also shows that $\mathcal M_{2d}$ is smooth of relative dimension 19.
\end{proof}
\begin{rem}
Since smoothness will be essential for all our further considerations, unless explicitly stated, by $\Fk_{2d}$ (respectively $\mathcal M_{2d}$) we will mean the smooth stack $\Fk_{2d}\otimes_\Z \Z[\frac{1}{2d}]$ (respectively $\mathcal M_{2d}\otimes_\Z \Z[\frac{1}{2d}]$) over $\Spec(\Z[\frac{1}{2d}])$.
\end{rem}
We will end this section speculating about other possible moduli spaces and functors of polarized K3 surfaces. Note first that one could have started with a moduli functor $\Fk'_{2d}$ of (primitively) polarized K3 schemes of degree $2d$. The problem we came up with restricting only to schemes was proving effectiveness of descent for K3 schemes. For this reason one takes the ``\'etale sheafification'' of $\Fk'_{2d}$ considering (primitively) polarized K3 spaces. This makes the descent obstruction essentially trivial.

Next, one can consider deformations of polarized K3 surfaces as in Section \ref{DefTh} by algebraic spaces and not only schemes. For a polarized K3 surface $(X_0, \l_0)$ over an algebraically closed field $k$ define
\begin{displaymath}
 \DefF_{\rm AlgSp}(X_0,\l_0) \colon \underline A \lr {\rm Sets}
\end{displaymath}
to be the functor sending an object $A$ of $\underline A$ to the isomorphism classes of triples $(\X,\l, \phi_0)$ where $(\X \lr \Spec(A), \l)$ is a polarized K3 space and $\phi_0$ is an isomorphism $\phi_0 \colon (\X,\L)\otimes_A k \cong (X_0,\L_0)$. Combining Theorem \ref{F-Prp}, Lemma \ref{smoothstac} and \cite[Cor. 10.11]{L-MB} we conclude that $\DefF_{\rm AlgSp}$ is pro-representable, formally smooth and of dimension $19$.
\section{Level Structures of Polarized K3 Surfaces}
Recall that for an abelian scheme $(A,\l)$ over a base scheme $S$ and a natural number $n$ which is invertible in $S$ one defines a (Jacobi) level $n$-structure on $A$ to be an isomorphism $\theta \colon A[n] \lr (\Z/n\Z)_S$ of \'etale sheaves on $S$ satisfying some further properties. In other words, one uses the Tate module of an abelian variety in order to define level structures. For a K3 surface $X$ we will use the same idea applied to $H^2_{\rm et}(X_{\bar k},\Z_l(1))$. More precisely, we will introduce the notion of level structures on primitively polarized K3
surfaces of degree $2d$ corresponding to open compact subgroups of $\SO(V_{2d}, \psi_{2d})(\A_f)$ (see below) and define moduli spaces of primitively polarized K3 surfaces with level structures. We set up some notations first.
\begin{itemize}
\item All schemes in this section will be assumed to be locally noetherian. 
\item For a finite set of primes $\mathsc B = \{p_1,\dots,p_r\}$ we denote by $\Z_{\mathsc B}$ the product $\prod_{p \in \mathsc B} \Z_p$ and by $N_{\mathsc B}$ the product of the primes in $\mathsc B$. 
\item We fix a natural number $d$. We shall use the notations $L_{2d,\mathsc B}$ and $L_{0,\mathsc B}$ for the quadratic lattices $L_{2d} \otimes \Z_{\mathsc B}$ and $L_0 \otimes \Z_{\mathsc B}$ (cf. Section \ref{QLatK3}).
\item Let $\Kg \subset \SO(V_{2d})(\hat \Z)$ be a subgroup of finite index and let $\mathsc B = \{p_1,\dots,p_r\}$ be the set of prime divisors of $2d$ and primes $p$ for which $\Kg_p \ne \SO(V_{2d})(\Z_p)$. We denote by $\Kg_{\mathsc B}$ the product $\prod_{p \in \mathsc B} \Kg_p$.
\end{itemize}
\noindent
\subsection{Level Structures}\label{LevelStructures}
Let $S$ be a connected scheme over $\Z[\frac{1}{p_1\dots p_r}]$
and suppose given a polarized K3 space $(\pi \colon X \lr
S, \l)$ of degree $2d$. Let $P^2_{\rm et}\pi_*\Z_{\mathsc B}(1)$ be
the sheaf of primitive cohomology i.e., the orthogonal complement
of $c_1(\l)$ in $R^2_{\rm et}\pi_*\Z_{\mathsc B}(1)$. Take a geometric point $\bar b$ of $S$
and let $\bar b \colon \Spec(k(\bar b)) \lr S$ be the
corresponding morphism of schemes. Consider the free $\Z_{\mathsc
B}$-module of rank 21
\begin{displaymath}
 P^2(\bar b) := \bar b^*P^2_{\rm et}\pi_*\Z_{\mathsc B}(1)
\end{displaymath}
i.e., the fiber of $P^2_{\rm et}\pi_*\Z_{\mathsc B}(1)$ at $\bar b$
with its action of $\pi_1^{\rm alg}(S, \bar b)$ and the bilinear form
$\psi_{\l,\Z_{\mathsc B}}$.

Suppose given an class $\a_{\bar b}$ in the set
\begin{displaymath}
 \biggl\{\Kg_{\mathsc B} \backslash \textrm{Isometry}\bigl(L_{2d, \Z_{\mathsc B}}, P^2(\bar
 b)\bigr)\biggr\}^{\pi_1^{\rm alg}(S,\bar b)}
\end{displaymath}
where $\Kg_{\mathsc B}$ acts on $\textrm{Isometry}\bigl(L_{2d, \Z_{\mathsc B}},
P^2(\bar b)\bigr)$ on the right via its action on $L_{2d, \Z_{\mathsc B}}$ and $\pi_1^{\rm alg}(S,\bar b)$ acts on the left via its action on $P^2(\bar b)$. Let $\bar b'$ be another geometric point in $S$. The $\a_{\bar b}$ determines uniquely a class in 
\begin{displaymath}
\biggl\{\Kg_{\mathsc B} \backslash \textrm{Isometry}\bigl(L_{2d, \Z_{\mathsc B}}, P^2(\bar b')\bigr)\biggr\}^{\pi_1^{\rm alg}(S,\bar b')} 
\end{displaymath}
in the following way: One can find an isomorphism 
\begin{equation}\label{etfundgr}
\delta_\pi \colon \pi_1^{\rm alg}(S,\bar b) \cong \pi_1^{\rm alg}(S,\bar b')
\end{equation}
and an isometry 
$$
\delta_{\rm et} \colon H^2_{\rm et}(X_{\bar b}, \Z_{\mathsc B}(1)) \lr H^2_{\rm et}(X_{\bar b'}, \Z_{\mathsc B}(1))
$$ 
determined uniquely by $\delta_\pi$, mapping $c_1(\l_{\bar b})$ to $c_1(\l_{\bar b'})$, such that $\delta_{\rm et}(\gamma\cdot x) = \delta_\pi(\gamma)\cdot \delta_{\rm et}(x)$ for every $x \in H^2_{\rm et}(X_{\bar b}, \Z_{\mathsc B}(1))$ and $\gamma \in \pi_1^{\rm alg}(S,\bar b)$. The isometry $\delta_{\rm et}$ defines an isometry between $P^2(\bar b)$ and $P^2(\bar b')$ which we will denote again by $\delta_{\rm et}$. Let $\tilde \a$ be a representative of the class $\a_{\bar b}$. Then the class $\a_{\bar b'}$ of $\delta_{\rm et} \circ \tilde \a$ in $\Kg_{\mathsc B} \backslash \textrm{Isometry}\bigl(L_{2d,\mathsc B}, P^2(\bar b')\bigr)$ is $\pi_1^{\rm alg}(S, \bar b')$-invariant. Any other representative $\tilde \a_1$ of $\a_{\bar b}$ differs by an element in $\Kg_{\mathsc B}$ and hence gives rise to the same class $\a_{\bar b'}$ in $\Kg_{\mathsc B} \backslash \textrm{Isometry}\bigl(L_{2d,\mathsc B}, P^2(\bar b')\bigr)$.

Any two isomorphisms \eqref{etfundgr} differ by an inner automorphism of $\pi_1^{\rm alg}(S,\bar b)$ and therefore we see that that class of $\delta_{\rm et} \circ \tilde \a$ is independent of the choice of an isomorphism \eqref{etfundgr}. 

This remark allows us to make the following definition.
\begin{dfn}\label{DefLevStr}
A \emph{level $\Kg$-structure} on a primitively polarized K3 space $(\pi \colon X \lr S, \l)$ over a connected scheme $S \in (\Sch/\Z[1/p_1 \dots p_r])$ is an element of the set
$$
 \biggl\{\Kg_{\mathsc B} \backslash \textrm{Isometry}\bigl(L_{2d, \Z_{\mathsc B}}, P^2(\bar
 b)\bigr)\biggr\}^{\pi_1^{\rm alg}(S,\bar b)}.
$$
The group $\Kg_{\mathsc B}$ acts on $\textrm{Isometry}\bigl(L_{2d, \Z_{\mathsc B}},
P^2(\bar b)\bigr)$ on the right via its action on $L_{2d, \Z_{\mathsc B}}$ and
 $\pi_1^{\rm alg}(S,\bar b)$ acts on the left via its action on $P^2(\bar b)$. In general,
 a level $\Kg$-structure on $(\pi \colon X \lr S, \l)$ is a level $\Kg$-structure on each
 connected component of $S$.
\end{dfn}
If $\tilde \a \colon L_{2d, \mathsc B} \lr P^2_{et}(\bar b)$ is a representative of the class $\a$, then via the isomorphism 
$$
 \tilde \a^{\rm ad} \colon {\rm O}(V_{2d})(\Z_{\mathsc B}) \cong {\rm O}(P^2(\bar b))(\Z_{\mathsc B})$$ 
the monodromy action 
$$
 \rho \colon \pi_1^{\rm alg}(S,\bar b) \lr {\rm O}(P^2(\bar b))(\Z_{\mathsc B})
$$ 
factorizes through ${\tilde \a}^{\rm ad}(\Kg_{\mathsc B})$.
\begin{rem}
If all residue fields of the points in $S$ in Definition \ref{DefLevStr} are of characteristic zero, then one can define a level $\Kg$-structure to be an element of set 
$$
 \biggl\{\Kg \backslash \textrm{Isometry}\bigl(L_{2d, \hat \Z}, P^2(\bar b)\bigr)\biggr\}^{\pi_1^{\rm alg}(S,\bar b)}
$$ 
where $P^2(\bar b) := \bar b^*P^2_{\rm et}\pi_*\hat \Z(1)$.
\end{rem}
We will consider two important examples of level structures on primitively polarized K3 spaces.
\begin{exa}\label{LnStr}
Fix a natural number $n$ and consider the group
\begin{displaymath}
 \Kg_n = \bigl\{ \gamma \in \SO(V_{2d})(\hat \Z) | \ \gamma
 \equiv 1 \pmod n \bigr\}.
\end{displaymath}
Then the set $\mathsc B$ consists of the prime divisors of $2dn$. We will give a direct interpretation of level $\Kg_n$-structures.

Let $S$ be a scheme over $\Z[1/2dn]$ and consider a primitively polarized K3 space $(\pi \colon X \lr S, \l)$ of degree $2d$.
As usual we denote by $P^2_{\rm et}\pi_*(\Z/n\Z)(1)$ the orthogonal
complement of $c_1(\l)$ in $R^2_{\rm et}\pi_*(\Z/n\Z)(1)$ with respect to the bilinear form $\psi_n = \psi \otimes_\Z
\Z/n\Z$. Then a level $\Kg_n$-structure amounts to giving an isomorphism
$$
 \a_N \colon \bigl(P^2_{\rm et}\pi_*(\Z/n\Z)(1), \psi_{\L, n}\bigr) \lr (L_{2d, \Z/n\Z},\psi_{2d, \Z/n\Z})_S
$$
of \'etale sheaves on $S$, where $(L_{2d, \Z/n\Z},\psi_{2d, \Z/n\Z})_S$ is the constant
polarized \'etale sheaf over $S$ with fibers $(L_{2d},\psi_{2d})\otimes \Z/n\Z$.

We will call level a $\Kg_n$-structure on $X$ simply a \emph{level
$n$-structure}.
\end{exa}
\begin{exa}\label{LnSpStr}
Let $G$ be the algebraic group $\SO(V_{2d})$ over $\Q$. Consider the \emph{even Clifford algebra}
$C^+(V_{2d}, \psi_{2d})$ over $\Q$ and let $G_1$ be the \emph{even Clifford
group} over $\Q$. In other words we set
\begin{displaymath}
 G_1 = \CSpin(V_{2d}) = \bigl\{g\in C^+(V_{2d})^*\ | \ g V_{2d} g^{-1} = V_{2d}\bigr\}.
\end{displaymath}
The natural homomorphism of linear algebraic groups $G_1 \lr
G$ given by $g \mapsto (v \mapsto gvg^{-1})$ fits into an exact sequence (see \cite[\S3.2]{D-K3})
$$
 0 \lr \mathbb G_m \lr G_1 \lr G \lr 0.
$$
Set $G_1(\Z)$ to be $G_1(\Q) \cap C^+(L_{2d})^*$. We have an exact sequence
(see \cite[\S 4.4]{A-HV})
\begin{equation}\label{IntCSpinExSq}
 0 \lr \Z/2\Z \lr G_1(\Z) \lr G(\Z).
\end{equation}
For a natural number $n$ denote
\begin{displaymath}
 \Gamma_n = \bigl\{ \gamma \in G(\Z)\ | \ \gamma \equiv 1 \pmod n \bigr\}
\end{displaymath}
and
\begin{displaymath}
 \Gamma_n^{\rm sp} = \bigl\{ \gamma \in G_1(\Z)\ |\ \gamma \equiv 1 \pmod n \ {\rm in} \ C^+(L_{2d}) \bigr\}.
\end{displaymath}
If $n > 2$, then $\Gamma_n$ and $\Gamma^{\rm sp}_n$ are torsion free. Hence one sees
from the exact sequence (\ref{IntCSpinExSq}) that $\Gamma_n^{\rm sp}$ is
isomorphic with its image $\Gamma_n^{\rm a}$ in $G(\Q)$.

Consider the group
\begin{displaymath}
 \Kg^{\rm sp}_n = \bigl\{\gamma \in G_1(\hat \Z)\ | \gamma \equiv 1 \pmod n\ {\rm in}\
 C^+(L_{2d, \hat \Z})\bigr\}.
\end{displaymath}
We have that $\Kg^{\rm sp}_n \cap G_1(\Q) = \Gamma_n^{\rm sp}$. Moreover the image
$\Kg^a_n$ of $\Kg^{\rm sp}_n$ in $G(\hat \Z)$ is of finite index. Indeed, for
every $l$ not dividing $2nd$, the $l$-component of $\Kg^{\rm a}_n$ is $G(\Z_l)$
as shown in \cite[\S 4.4]{A-HV}. Hence the set $\mathsc B$ for $\Kg^{\rm a}_n$ is the set of prime divisors of $2dn$.

We consider polarized K3 surfaces with level $\Kg^{\rm a}_n$-structure. Note that this
level structure is in general finer than level $\Kg_n$-structure as $\Kg^{\rm a}_n \subset \Kg_n$ is of finite index. We will call it \emph{spin level $n$-structure}.
\end{exa}

\subsection{Motivation}
We will pause here and give a motivation for the rest of the definitions we make in this section. So far we have defined level $\Kg$-structures using the primitive second \'etale cohomology group of a polarized K3 surface. Using these level structures one can define moduli stacks $\Fk_{2d,\Kg}$ of primitively polarized K3 surfaces of degree $2d$ with a level $\Kg$-structure and show that they are algebraic spaces (cf. Theorem \ref{MStLevStr} below). Over $\C$, we can relate these spaces to the orthogonal Shimura variety associated to the group $\SO(2,19)$. More precisely in Chapter 3, Section 3.4.2 of \cite{JR-Thesis} we define a period morphism
$$
 j_{d,\Kg,\C} \colon \Fk_{2d,\Kg,\C} \lr Sh_{\Kg}(\SO(2,19),\Omega^\pm)_\C
$$
which is \'etale. This is similar to the case of moduli of abelian varieties where one can identify $\Av_{g,1,n}\otimes \C$ with $Sh_{\Lambda_n}(\CSp_{2g},\mathfrak H_{g}^\pm)_\C$. In general, due to the fact that the injective homomorphism \eqref{InjSOHom}
$$
 i^{\rm ad} \colon \{g \in \SO(V_0)(\Z)\ |\ g(e_1-df_1) = e_1-df_1\} \hookrightarrow \SO(V_{2d})(\Z)
$$ 
defined in Section \ref{QLatK3} is not surjective, the period map $j_{d,\Kg,\C}$ need not be injective. In order to construct an injective period morphism we will define level structures using the ``full'' second \'etale cohomology group of a K3 surface. We use these full level structures in \cite[Ch. 3]{JR-Thesis} to show that every complex K3 surface with complex multiplication by a CM-field $E$ is defined over an abelian extension of $E$.
\newline
\newline 
\subsection{Full Level Structures} 
The inclusion of lattices $i \colon L_{2d} \hookrightarrow L_0$ (see Section \ref{QLatK3}) defines injective homomorphisms of groups
\begin{displaymath}
i^{\rm ad} \colon \{g \in {\rm O}(V_0)(\hat \Z)\ |\ g(e_1-df_1) = e_1-df_1\} \hookrightarrow {\rm O}(V_{2d})(\hat \Z)
\end{displaymath}
and
\begin{displaymath}
i^{\rm ad} \colon \{g \in \SO(V_0)(\hat \Z)\ |\ g(e_1-df_1) = e_1-df_1\} \hookrightarrow \SO(V_{2d})(\hat \Z).
\end{displaymath}

\begin{dfn}\label{AdmissGrDef}
A subgroup $\Kg \subset \SO(V_{2d})(\hat \Z)$ of finite index is called \emph{admissible} if it is contained in the image 
\begin{displaymath}
i^{\rm ad}\bigl(\{g \in \SO(V_0)(\hat \Z)\ |\ g(e_1-df_1) = e_1-df_1 \}\bigr) \subset \SO(V_{2d})(\hat \Z).
\end{displaymath}
\end{dfn}
If $\Kg$ is an admissible subgroup of $\SO(V_{2d})(\hat \Z)$ then all its subgroups of finite index $\Kg'\subset \Kg$ are also admissible. 
\begin{exa}
The group $\Kg_{2d}$ is admissible. Hence all its subgroups of finite index are admissible, as well.
\end{exa}
\begin{exa}
If $d = 1$ then $\Kg_n$ is admissible for any $n \geq 2$.
\end{exa}
Let $\Kg$ be an admissible subgroup of $\SO(V_{2d})(\hat \Z)$ and let $\mathsc B$ be the set, consisting of all prime divisors of $2d$ and, of the primes $p$ for which $\Kg_p \ne \SO(V_{2d})(\Z_p)$. Using the notations introduced before Definition \ref{DefLevStr} we set
\begin{displaymath}
 H^2(\bar b) :=  \bar b^*R^2_{\rm et}\pi_*\Z_{\mathsc B}(1).
\end{displaymath}
In order to simplify the notations we will identify a subgroup of $\{g \in \SO(V_0)(\hat \Z)\ |\ g(e_1-df_1) = e_1-df_1\}$ with its image in $\SO(V_{2d})(\hat \Z)$ under the injective homomorphism $i^{\rm ad}$.
\begin{dfn}
A \emph{full level $\Kg$-structure} on a primitively polarized K3 space $(\pi \colon X \lr S, \l)$ over a connected scheme $S \in (\Sch/\Z[1/p_1 \dots p_r])$ is an element of the set
$$
 \biggl\{\Kg_{\mathsc B} \backslash \bigl\{g\in \textrm{Isometry}\bigl(L_{0,\Z_{\mathsc B}}, H^2(\bar b)\bigr)\ |\ g(e_1-df_1) = c_1(\l_{\bar b})\bigr\}\biggr\}^{\pi_1^{\rm alg}(S,\bar b)}.
$$
The group $\Kg_{\mathsc B}$ acts on $\bigl\{g \in \textrm{Isometry}\bigl(L_{0, \Z_{\mathsc B}},
H^2(\bar b)\bigr)\ |\ g(e_1-df_1) = c_1(\l_{\bar b})\bigr\}$ on the right via its action on $L_{0, \Z_{\mathsc B}}$ and $\pi_1^{\rm alg}(S,\bar b)$ acts on the left via its action on $H^2(\bar b)$. A full level $\Kg$-structure on $(\pi \colon X \lr S, \l)$ over a general base $S$ is a full level $\Kg$-structure on each connected component of $S$.
\end{dfn}
Again, a class $\a_{\bar b}$ for a geometric point $\bar b$ as above determines uniquely a class $\a_{\bar b'}$ for any other geometric point $\bar b'$. If $\tilde \a \colon L_{0, \mathsc B} \lr H^2(\bar b)$ is a representative of the class $\a$, then via the isomorphism 
$$
 \tilde \a^{\rm ad} \colon  {\rm O}(V_{0})(\Z_{\mathsc B}) \cong {\rm O}(H^2(\bar b))(\Z_{\mathsc B})
$$
the monodromy action $\rho \colon \pi_1^{\rm alg}(S,\bar b) \lr {\rm O}(H^2(\bar b))(\Z_{\mathsc B})$ factorizes through ${\tilde \a}^{\rm ad}(\Kg_{\mathsc B})$.
\begin{exa}\label{FullLNStr}
Let $n\geq 3$ be an integer. Define the group
\begin{displaymath}
\Kg_n^{\rm full} = \bigl\{g \in \SO(V_0)(\hat \Z) |\ g(e_1 -df_1) = e_1-df_1\ {\rm and}\ g
 \equiv 1 \pmod n \bigr\}.
\end{displaymath}
By definition it is an admissible subgroup of $\SO(V_{2d})(\hat \Z)$. Let $S$ be a scheme over $\Z[1/2dn]$ and consider a K3 space $(\pi \colon X \lr S, \l)$ with a primitive polarization of degree $2d$. Then a full level $\Kg_n^{\rm full}$-structure amounts to giving an isomorphism
$$
 \a_N \colon \bigl(R^2_{\rm et}\pi_*(\Z/n\Z)(1), \psi\bigr) \lr (L_{0, \Z/n\Z},\psi_{0, \Z/n\Z})_S
$$
of \'etale sheaves on $S$, where $(L_{0, \Z/n\Z},\psi_{0, \Z/n\Z})_S$ is the constant
polarized \'etale sheaf over $S$ with fibers $(L_0,\psi_0)\otimes \Z/n\Z$.

We will call a full level $\Kg_n^{\rm full}$-structure on $X$ simply a \emph{full level $n$-structure}.
\end{exa}
\section{Moduli Spaces of Polarized K3 Surfaces with a Level Structure}\label{MSPK3LS}
In this section we will use the notion of a (full) level structure level structure to define moduli functors of primitively polarized K3 spaces with a (full) level structure. Using Artin's criterion and Proposition \ref{finaut} we will show that these functors are representable by algebraic spaces over open parts of $\Spec(\Z)$.

We shall be using the notations established in the beginning of Section \ref{LevelStructures}. In particular we fix a natural number $d$. To a subgroup $\Kg$ of $\SO(V_{2d})(\hat \Z)$ we associated a finite set of primes $\mathsc B$ and $N_{\mathsc B}$ will denote the product of these primes.
\subsection{Moduli of K3 Surfaces with Level Structure} 
Let $\Kg$ be a subgroup of $\SO(\hat \Z)$ of finite index. We will assume further that it is contained in $\Kg_n$ for some $n \geq 3$. Let $\X_1 = (\pi_1 \colon X_1 \lr S_1, \l_1)$ and $\X_2 = (\pi_2 \colon X_2 \lr S_2)$ be two objects of $\Fk_{2d}$. Suppose that $S_1$ and $S_2$ are connected and let $(f,f_S) \in \Hm(\X_1, \X_2)$ (in $\Fk_{2d}$). Let $\bar b_1$ and $\bar b_2$ be two geometric points of
$S_1$ and $S_2$ such that $f_S(\bar b_1) = \bar b_2$. Then the
morphism $f$ defines a homomorphism $f^*_{\rm et} \colon P^2(\bar b_2) \lr P^2(\bar b_1)$. Hence we obtain a map
\begin{displaymath}
 f^\vee \colon \Kg_{\mathsc B} \backslash \textrm{Isometry} \bigl(L_{2d,\Z_{0,\mathsc B}}, P^2(\bar
 b_2)\bigr) \lr \Kg_{\mathsc B} \backslash \textrm{Isometry} \bigl(L_{2d,\Z_{\mathsc B}}, P^2(\bar
 b_1)\bigr)
\end{displaymath}
given by $\a \mapsto f^*_{\rm et} \circ \a$ and commuting with the monodromy actions on both sides.
\begin{dfn}
For $d$ and $\Kg$ as above consider the category $\Fk_{2d, \Kg}$ defined in the following way:
\begin{enumerate}
\item [{\it Ob:}] Triples $(\pi \colon X \lr S, \l, \a)$ of a K3 space $\pi \colon X \lr S$ with a primitive polarization $\l$ of degree $2d$ and with a level $\Kg$-structure $\a$ on $(\pi \colon X \lr S, \l)$.
\item [{\it Mor:}] Suppose given two triples $\X_1 = (\pi_1 \colon X_1 \lr S_1, \l_1,
\a_1)$ and $\X_2 = (\pi_2 \colon X_2 \lr S_2, \l_2, \a_2)$.  Let $f_S \colon S_1 \lr S_2$
be a morphism of schemes. Choose base geometric points $\bar
b'_1$ and $\bar b'_2$ on any two connected components $S'_1$ and $S'_2$ of $S_1$
and $S_2$ for which $f \colon S'_1 \lr S'_2$ such that
$f_S(\bar b'_1) = \bar b'_2$. Define the morphisms between $\X_1$ and $\X_2$ in the following way
\begin{displaymath} 
\begin{split}
 \Hm(\X_1, \X_2) = \bigl\{{\rm pairs}\ (f_S,f)\ |& f_S \colon S_1 \lr S_2 \ {\rm is\ a\
 morph.\ of\ spaces,} \\
  & f \colon X_1 \lr X_2\times_{S_2,f_S} S_1\ {\rm is\ an\ isom.\ of} \\
  & S_1-{\rm spaces}\ {\rm with}\ f^*\l_2 = \l_1\ {\rm and}\\
  & f^\vee (\a_1) = \a_2\ {\rm on\ any\ conn.\ cmpt.\ of}\ S_1 \bigr\} \\
\end{split}
\end{displaymath}
\end{enumerate}
\end{dfn}
\noindent
Next we define three projection functors.
\begin{enumerate}
\item [1.] Consider the following forgetful functor
\begin{displaymath}
 pr_{\Fk_{2d,\Kg}} \colon \Fk_{2d, \Kg} \lr (\Sch/\Z[1/N_{\mathsc B}])
\end{displaymath}
sending a triple $(\pi \colon X \lr S, \l, \a)$ to $S$. It makes $\Fk_{2d, \Kg}$ into a category over $(\Sch/\Z[1/N_{\mathsc B}])$.
\item [2.] For any $\Kg$, satisfying the assumptions of the beginning of the section, one has a projection functor
\begin{equation}\label{prK}
 pr_\Kg \colon \Fk_{2d,\Kg} \lr \Fk_{2d, \Z[1/N_{\mathsc B}]}
\end{equation}
sending a triple $(\pi \colon X \lr S, \l, \a)$ to $(\pi \colon X \lr S, \l)$
and an element $(f,f_S)\in \Hm(\X,\mathcal Y)$ of $\Fk_{2d,\Kg}$ to $(f,f_S)$.
\item [3.] For any two subgroups $\Kg_1 \subset \Kg_2$ of finite index in $\SO(V_{2d})(\hat \Z)$ (contained in some $\Kg_n$ for $n \geq 3$) one has a projection functor
\begin{equation}\label{prK1K2}
 pr_{(\Kg_1,\Kg_2)} \colon \Fk_{2d,\Kg_1,\Z[1/N_{\mathsc B_1 \cup \mathsc B_2}]} \lr \Fk_{2d,\Kg_2,\Z[1/N_{\mathsc B_1 \cup \mathsc B_2}]}.
\end{equation}
It sends an object $(X\lr S, \l,\a_{\Kg_1})$ to $(X\lr S, \l,\a_{\Kg_2})$ where $\a_{\Kg_2}$ is the class of $\a_{\Kg_1}$ in $\Kg_{2, \mathsc B} \backslash \textrm{Isometry}\bigl(L_{2d, \Z_{\mathsc B}}, P^2(\bar b)\bigr)$. Morphism of $\Fk_{2d,\Kg_1,\Z[1/N_{\mathsc B_1 \cup \mathsc B_2}]}$ are mapped to morphism of $\Fk_{2d,\Kg_2,\Z[1/N_{\mathsc B_1 \cup \mathsc B_2}]}$ in the obvious way.
\end{enumerate}
From the definitions of the functors we see that $pr_{\Kg_1} = pr_{(\Kg_1,\Kg_2)} \circ pr_{\Kg_2}$ over $\Z[1/N_{\mathsc B_1 \cup \mathsc B_2}]$. 
\begin{thm}\label{MStLevStr}
The category $\Fk_{2d,\Kg}$ is a separated algebraic space over $\Z[1/N_{\mathsc B}]$. It is smooth of relative dimension 19 and the forgetful morphism \eqref{prK}
$$
 pr_\Kg \colon \Fk_{2d,\Kg} \lr \Fk_{2d, \Z[1/N_{\mathsc B}]}
$$
is finite and \'etale.
\end{thm}
\begin{proof} We divide the proof into several steps.
\newline
{\bf Step 1:} The category $\Fk_{2d,\Kg}$ is a stack. The proof goes exactly in the same lines
as the one of Lemma \ref{LStack}. We will use Artin's criterion (cf. \cite[Cor. 10.11]{L-MB}) to show that $\Fk_{2d,\Kg}$ is an algebraic space.

We claim that the diagonal morphism $\Delta \colon \Fk_{2d,\Kg} \lr \Fk_{2d,\Kg} \times_{\Z[1/N_{\mathsc B}]} \Fk_{2d,\Kg}$ is representable, separated and of finite type. By Remark 4.1.2 in \cite{L-MB} it is equivalent to showing that for any two objects $\X=(X\lr S,\l_X,\a_X)$ and $\mathcal Y = (Y \lr S,\l_Y,\a_Y)$ the functor $\Is_S(\X,\mathcal Y)$ has these properties. We will prove first the following result.
\begin{lem}
For any object $\X$ of $\Fk_{2d,\Kg}$ we have that $\Aut_S(\X) = \{\id_\X\}$.
\end{lem}
\begin{proof}
By assumption the group $\Kg$ is contained in $\Kg_n$ for some $n \geq 3$. Hence a level $\Kg$-structure on a primitively polarized K3 space $(X \lr S, \l)$ defines in a natural way (using the functor $pr_{(\Kg,\Kg_n)}$) a level $n$-structure $\a_n$ on $X$. We have that 
$$
 \Aut_S\bigl((X \lr S, \l, \a)\bigr)(U) \subset \Aut_S\bigl((X \lr S, \l, \a_n)\bigr)(U)
$$
for an $S$-scheme $U$ hence it is enough to prove the lemma assuming that $\Kg = \Kg_n$. 
 
Let $\X = (X \lr S, \l, \a)$ be an object in $\Fk_{2d,\Kg}$, let $f \in \Aut_S(\X)(U)$ and assume that $U$ is connected. Take a geometric point $\bar b \colon \Spec(\Omega) \lr U$. Then for the finite set $\mathsc B = \{{\rm the\ prime\ divisors\ of}\ $n$\}$ the morphism $f$ induces an automorphism 
$$
 f^*_{\rm et} \colon H^2_{\rm et}(X_{\bar b}, \Z_{\mathsc B}(1)) \lr H^2_{\rm et}(X_{\bar b}, \Z_{\mathsc B}(1))
$$ 
fixing $c_1(\l_{\bar b})$ and such that
$$
 f^*_{\rm et} \colon P^2_{\rm et}(X_{\bar b}, \Z/n\Z(1)) \lr P^2_{\rm et}(X_{\bar b}, \Z/n\Z(1))
$$
is the identity (cf. Example \ref{LnStr}). As the automorphism $f$ is of finite order we have that $f^*_{\rm et} \in {\rm O}\bigl(P^2_{\rm et}(X_{\bar b},\Z_{\mathsc B}(1))\bigr)$ is semi-simple and its eigenvalues are roots of unity. We have further that $f^*_{\rm et} \equiv 1 \pmod n$ so we conclude by \cite[Ch. IV, Application II, p. 207, Lemma]{Mum-AV} that $f^*_{\rm et}$ is the identity automorphism of $P^2_{\rm et}(X_{\bar b},\Z_{\mathsc B}(1))$. As it fixes $c_1(\l_{\bar b})$ we see that it acts as the identity on $H^2_{\rm et}(X_{\bar b},\Z_{\mathsc B}(1))$. Therefore by Proposition \ref{finaut} we that $f = \id_{X_{\bar b}}$. As the geometric point $\bar b$ can be chosen arbitrary we have that $f = \id_{\X_U}$.
\end{proof}
We see from the lemma that for a $S$-scheme $U$ the set $\Is_S(\X,\mathcal Y)(U)$ is either empty or it consists of one element. Indeed, suppose that $f_i \in \Is_S(\X,\mathcal Y)(U)$ for $i =1,2$. Then the composition $f^{-1}_2\circ f_1$ belongs to $\Aut_S(\X)(U)$ and hence it is the identity. This shows that $\Is_S(\X,\mathcal Y)$ is representable and of finite type. The fact that it is unramified and separated over $S$ follows from Lemma \ref{Isom} as one has that 
$$
\Is_S(\X,\mathcal Y)(U) \subset \Is_S\bigl((X\lr S,\l_X),(Y\lr S,\l_Y)\bigr)(U).
$$ 
Next we claim that the stack $\Fk_{2d,\Kg}$ is locally of finite presentation. This follows from \cite[Expos\'e IX, 2.7.4]{SGA3} and the fact that $\Fk_{2d}$ is locally of finite presentation. Conditions {\it (iii)} and {\it (iv)} of \cite[Cor. 10.11]{L-MB} follow from the corresponding properties of $\Fk_{2d}$ and the fact that for any small surjection of rings $R \lr R'$ the category of \'etale schemes over $R$ is equivalent to the category of \'etale schemes over $R'$ (\cite[EGA IV, 18.1.2]{EGA}).

Thus $\Fk_{2d,\Kg}$ is an algebraic stack. As $\Aut_S(\X) = \{\id_\X\}$ for any object we have that $\Fk_{2d,\Kg}$ is an algebraic space (\cite[Cor. 8.1.1]{L-MB}).
\newline
\newline
{\bf Step 2:} We will show that the morphism of algebraic stacks $pr_\Kg \colon \Fk_{2d,\Kg}
\lr \Fk_{2d, \Z[1/N_{\mathsc B}]}$ is representable and \'etale. Indeed, let $S$ be
a connected scheme and suppose given a morphism $S \lr \Fk_{2d}$ i.e., a
polarized K3 space $(\pi \colon X \lr S, \l)$ over $S$. Let $\bar b \colon \Spec(\Omega) \lr S$ be a geometric point of $S$. Let $\rho \colon \pi^{\rm alg}(S,\bar b) \lr {\rm O}(P^2(\bar b))$ be the monodromy representation and let $\tilde a \colon L_{2d,\mathsc B} \lr P^2(\bar b)$ be an isometry. Then the preimage $\rho^{-1} \circ \a^{\rm ad} (\Kg_{\mathsc B})$ is an open subgroup of $\pi^{\rm alg}_1(S,\bar b)$ (of finite index) and hence it defines an \'etale cover $S_{\tilde \a}$ of $S$. One has that the class $\a$ of $\tilde \a$ in $\Kg_{\mathsc B} \backslash \textrm{Isometry} \bigl(L_{2d,\Z_{\mathsc B}}, P^2(\bar b)\bigr)$ is $\pi^{\rm alg}_1(S_{\tilde \a},\bar b)$-invariant by construction (for a fixed geometric point $\bar b \in S_{\tilde \a}$ over $\bar b$). Therefore we obtain a primitively polarized K3 space $(X_{S_{\tilde \a}} \lr S_{\tilde \a}, \l_{S_{\tilde \a}}, \a)$ with a level $\Kg$-structure $\a$. For two markings $\tilde \a_1$ and $\tilde \a_2$ we have that $\tilde \a_1^{\rm ad}(\Kg_{\mathsc B}) = \tilde \a_2^{\rm ad}(\Kg_{\mathsc B})$ if and only if $\tilde \a_2^{-1} \circ \tilde \a_1$ is an element of the normalizer $N_{{\rm O}(V_{2d})(\Z_{\mathsc B})}(\Kg_{\mathsc B})$ of $\Kg_{\mathsc B}$ in ${\rm O}(V_{2d})(\Z_{\mathsc B})$.

Denote by $S'$ the disjoint union of $S_{\tilde \a}$ where $\tilde \a$ runs over all (finitely many) classes in ${\rm O}(V_{2d})(\Z_{\mathsc B})/N_{{\rm O}(V_{2d})(\Z_{\mathsc B})}(\Kg_{\mathsc B})$. Let $(X' \lr S', \l_{S'}, \a)$ be the primitively polarized K3 space with a level $\Kg$-structure given by the triple $(X_{S_{\tilde \a}} \lr S_{\tilde \a}, \l_{S_{\tilde \a}}, \a)$ on the $\tilde \a$-{\it th} connected component $S_{\tilde \a}$ of $S'$. Then by construction we have a morphism of algebraic spaces 
$$
 \pi \colon S' \lr S\times_{\Fk_{2d,\Z[1/N_{\mathsc B}]}} \Fk_{2d,\Kg}
$$
over $S$. This morphism is surjective. Indeed, by \cite[Prop. 5.4]{L-MB} this condition can be checked on points, in which case it is obvious by construction. The morphism $S' \lr S$ is \'etale and therefore we conclude that $pr_{\Kg, S} \colon S\times_{\Fk_{2d,\Z[1/N_{\mathsc B}]}} \Fk_{2d,\Kg} \lr S$ and $\pi$ are also \'etale. Hence $pr_\Kg$ is \'etale. 
\newline
\newline
{\bf Step 3:} By Step 2 and Theorem \ref{F-Prp} the algebraic space $\Fk_{2d,\Kg}$ is smooth and of relative dimension 19 over $\Z[1/N_{\mathsc B}]$.
\end{proof}
\begin{rem}
Let $\Kg_1 \subset \Kg_2 \subset \Kg_n$ be subgroups of finite index in $\SO(V_{2d})(\hat \Z)$ and suppose that $n \geq 3$. Then the morphism \eqref{prK1K2} of algebraic spaces 
$$
pr_{(\Kg_1,\Kg_2)} \colon \Fk_{2d,\Kg_1,\Z[1/N_{\mathsc B_1 \cup \mathsc B_2}]} \lr \Fk_{2d,\Kg_2,\Z[1/N_{\mathsc B_1 \cup \mathsc B_2}]}
$$
is finite and \'etale. This follows from the theorem above and the relation $pr_{\Kg_1} = pr_{\Kg_2} \circ pr_{(\Kg_1,\Kg_2)}$.
\end{rem}
\begin{exa}
Let $n \geq 3$ be a natural number. Consider the group $\Kg_n$ defined in Example
\ref{LnStr}. We define $\Fk_{2d,n} = \Fk_{2d,\Kg_n}$ to be \emph{the moduli space of
primitively polarized K3 surfaces with level $n$-structure} over $\Z[1/2dn]$.
\end{exa}
\begin{exa}
Fix a natural number  $n \geq 3$  and consider the group $\Kg^{\rm a}_n$ defined in Example
\ref{LnSpStr}. We define $\Fk_{2d,n^{\rm sp}} = \Fk_{2d,\Kg^{\rm a}_n}$ to be \emph{the moduli space of
polarized K3 surfaces with spin level $n$-structure} over $\Z[1/2dn]$.
\end{exa}

\subsection{Moduli K3 Spaces with Full Level Structures}
Suppose that $\Kg \subset \Kg_n$ for some $n \geq 3$ is an admissible subgroup of $\SO(V_{2d})(\hat \Z)$. Let $\X_1 = (\pi_1 \colon X_1 \lr S_1, \l_1)$ and $\X_2 = (\pi_2 \colon X_2 \lr S_2)$ be two objects of $\Fk_{2d}$. Suppose that $S_1$ and $S_2$ are connected and let $(f,f_S) \in \Hm(\X_1, \X_2)$ (in $\Fk_{2d}$). Let $\bar b_1$ and $\bar b_2$ be two geometric points of
$S_1$ and $S_2$ such that $f_S(\bar b_1) = \bar b_2$. Then the
morphism $f$ defines a homomorphism $f^*_{\rm et} \colon H^2(\bar b_2) \lr H^2(\bar b_1)$ sending the class of $\l_{\bar b_2}$ to the class of $\l_{\bar b_1}$. Hence we obtain a map
\begin{gather*}
f^\vee \colon \Kg_{\mathsc B} \backslash \bigl\{\ g \in \textrm{Isometry} \bigl(L_{0,\Z_{\mathsc B}}, H^2(\bar b_2)\bigr)\ |\  g(e_1-df_1) = c_1(\l_{2, \bar b_2})\ \bigr\} \lr \\
\lr \Kg_{\mathsc B} \backslash \bigl\{\ g \in \textrm{Isometry} \bigl(L_{0,\Z_{\mathsc B}}, P^2(\bar b_1)\bigr)\ |\ g(e_1-df_1) = c_1(\l_{1,\bar b_1})\ \bigr\}
\end{gather*}
given by $\a \mapsto f^*_{et} \circ \a$ and commuting with the monodromy actions on both sides.
\begin{dfn}
For a natural number $d$ and an admissible subgroup $\Kg$ of $\SO(V_{2d})(\hat \Z)$ as above consider the category $\Fk_{2d, \Kg}^{\rm full}$ defined in the following way:
\begin{enumerate}
\item [{\it Ob:}] Triples $(\pi \colon X \lr S, \l, \a)$ of a K3 space $\pi \colon X \lr S$
over $S$ with a primitive polarization $\l$ of degree $2d$ and with a full level $\Kg$-structure $\a$ on $(\pi \colon X \lr S, \l)$.
\item [{\it Mor:}] Suppose given two triples $\X_1 = (\pi_1 \colon X_1 \lr S_1, \l_1,
\a_1)$ and $\X_2 = (\pi_2 \colon X_2 \lr S_2, \l_2, \a_2)$.  Let $f_S \colon S_1 \lr S_2$
be a morphism of schemes. Choose base geometric points $\bar
b'_1$ and $\bar b'_2$ on any two connected components $S'_1$ and $S'_2$ of $S_1$
and $S_2$ for which $f \colon S'_1 \lr S'_2$ such that
$f_S(\bar b'_1) = \bar b'_2$. Define the morphisms between $\X_1$ and $\X_2$ in the following way
\begin{displaymath} 
\begin{split}
 \Hm(\X_1, \X_2) = \bigl\{{\rm pairs}\ (f_S,f)\ |& f_S \colon S_1 \lr S_2 \ {\rm is\ a\
 morph.\ of\ spaces,} \\
  & f \colon X_1 \lr X_2\times_{S_2,f_S} S_1\ {\rm is\ an\ isom.\ of} \\
  & S_1-{\rm spaces}\ {\rm with}\ f^*\l_2 = \l_1\ {\rm and}\\
  & f^\vee (\a_1) = \a_2\ {\rm on\ any\ conn.\ cmpt.\ of}\ S_1 \bigr\} \\
\end{split}
\end{displaymath}
\end{enumerate}
\end{dfn}
\noindent
A full level $\Kg$-structure $\a$ on a primitively polarized K3 space $(X\lr S,\l)$ defines in a natural way a level $\Kg$-structure via the injective morphism
\begin{gather*}
i^\vee_{\Z_{\mathsc B}} \colon \Kg_{\mathsc B} \backslash \bigl\{\ g\in \textrm{Isometry}\bigl(L_{0,\Z_{\mathsc B}}, H^2(\bar b)\bigr)\ |\ g(e_1-df_1) = c_1(\l_{\bar b})\ \bigr\} \hookrightarrow \\
\hookrightarrow \Kg_{\mathsc B} \backslash \textrm{Isometry}\bigl(L_{2d, \Z_{\mathsc B}}, P^2(\bar b)\bigr)
\end{gather*}
commuting with the monodromy action. This morphism is defined by the embedding of lattices $i \colon L_{2d} \hookrightarrow L_0$ (see \eqref{InjOHom} in Section \ref{QLatK3}). Using this, just like in the case of moduli of primitively polarized K3 surfaces with a level structure, we define natural functors.
\begin{enumerate}
\item [4.] Define a functor
$$
 i_\Kg \colon \Fk_{2d,\Kg}^{\rm full} \lr \Fk_{2d,\Kg}   
$$
sending $(X\lr S,\l,\a)$ to $(X,\lr S,\l,i^\vee(\a))$ which makes $\Fk_{2d,\Kg}^{\rm full}$ into a full subcategory of $\Fk_{2d,\Kg}$ over $(\Sch/\Z[1/N_{\mathsc B}])$. 
\item [5.] One has the forgetful functor
\begin{equation}\label{prKfull}
 pr_\Kg \colon \Fk_{2d,\Kg}^{\rm full} \lr \Fk_{2d, \Z[1/N_{\mathsc B}]}
\end{equation}
sending a triple $(\pi \colon X \lr S, \l, \a)$ to $(\pi \colon X \lr S, \l)$
and an element $(f,f_S)\in \Hm(\X,\mathcal Y)$ of $\Fk_{2d,\Kg}^{\rm full}$ to $(f,f_S)$.
\item [6.] For any two admissible subgroups $\Kg_1 \subset \Kg_2$ of $\SO(V_{2d})(\hat \Z)$, contained in some $\Kg_n$ for $n \geq 3$, one has a projection functor
\begin{equation}\label{prK1K2full}
 pr_{(\Kg_1,\Kg_2)} \colon \Fk_{2d,\Kg_1,\Z[1/N_{\mathsc B_1 \cup \mathsc B_2}]}^{\rm full} \lr \Fk_{2d,\Kg_2,\Z[1/N_{\mathsc B_1 \cup \mathsc B_2}]}^{\rm full}
\end{equation}
defined in a similar way as the corresponding morphism \eqref{prK1K2} in 3.
\end{enumerate}
The functors $pr_\Kg$ and $pr_{(\Kg_1,\Kg_2)}$ defined above are the restrictions of the corresponding functors \eqref{prK} and \eqref{prK1K2} to the category of primitively polarized K3 surfaces with full level $\Kg_j$-structures via $i_{\Kg_j}$ for $j = 1,2$.
\begin{thm}\label{MSpaceFullLStr}
Let $\Kg$ be an admissible subgroup of $\SO(V_{2d})(\hat \Z)$ contained in $\Kg_n$ for some $n \geq 3$. The category $\Fk_{2d,\Kg}^{\rm full}$ is a separated, smooth algebraic space of relative dimension 19 over $\Z[1/N_{\mathsc B}]$. The morphism $p_{2d,\Kg} \colon \Fk_{2d,\Kg}^{\rm full} \lr \Fk_{2d,\Z[1/N_{\mathsc B}]}$ is \'etale and the morphism $i_\Kg \colon \Fk_{2d,\Kg}^{\rm full} \hookrightarrow \Fk_{2d,\Kg}$ is an open immersion.
\end{thm}
\begin{proof}
To prove that $\Fk_{2d,\Kg}^{\rm full}$ is representable by an algebraic space of finite type over $\Z[1/N_{\mathsc B}]$ one follows the steps of the proof of Theorem \ref{MStLevStr}. In this way we also see that the projection morphism $p_{2d,\Kg} \colon \Fk_{2d,\Kg}^{\rm full} \lr \Fk_{2d,\Z[1/N_{\mathsc B}]}$ is finite and \'etale. Therefore we have a commutative diagram
$$
\xymatrix{\Fk_{2d,\Kg}^{\rm full} \ar[dr]^{p_{2d,\Kg}} \ar[rr]^{i_\Kg} & & \Fk_{2d,\Kg} \ar[dl]^{p_{2d,\Kg}} \\
& \Fk_{2d,\Z[1/N_{\mathsc B}]} &
} 
$$
where the two morphisms $p_{2d,\Kg}$ are \'etale and surjective. Hence $i_\Kg$ is also \'etale and therefore it is open.
\end{proof}
\begin{rem}
Let $\Kg_1 \subset \Kg_2$ be two admissible subgroups of $\SO(V_{2d})(\hat \Z)$. Then the morphism of algebraic spaces 
$$
 pr_{(\Kg_1,\Kg_2)} \colon \Fk_{2d,\Kg_1,\Z[1/N_{\mathsc B_1 \cup \mathsc B_2}]}^{\rm full} \lr \Fk_{2d,\Kg_2,\Z[1/N_{\mathsc B_1 \cup \mathsc B_2}]}^{\rm full}
$$
is finite and \'etale. This follows from the theorem above and the relation $pr_{\Kg_1} = pr_{\Kg_2} \circ pr_{(\Kg_1,\Kg_2)}$.
\end{rem}
\begin{exa}
Let $n \geq 3$ be a natural number. Consider the group $\Kg_n^{\rm full}$ defined in Example
\ref{FullLNStr}. We define $\Fk_{2d,n}^{\rm full} = \Fk_{2d,\Kg_n^{\rm full}}^{\rm full}$ to be \emph{the moduli space of primitively polarized K3 surfaces with full level $n$-structure} over $\Z[1/2dn]$.
\end{exa}
\bibliographystyle{amsalpha}

\begin{thebibliography}{Mum74}

\bibitem[AGV71]{SGA3}
M.~Artin, A.~Grothendieck, and J-L. Verdier, \emph{{T}h\'eorie des {T}opos et
  {C}ohomologie {\'e}tale des {S}ch\'emas}, Lecture Notes in Mathematics, vol.
  269, 270, 305, Springer-Verlag, 1971.

\bibitem[And96]{A-HV}
Y.~Andr\'e, \emph{{O}n the {S}hafarevich and {T}ate {C}onjectures for
  {H}yperk\"ahler {V}arieties}, Mathematische Annalen \textbf{305} (1996),
  205--248.

\bibitem[Art69]{Art-AFM1}
M.~Artin, \emph{{A}lgebrization of {F}ormal {M}oduli: {I}}, in Global Analysis,
  Papers in honor of K. Kodaira, University of Tokyo Press/Princeton University
  Press, 1969, pp.~21--71.

\bibitem[BBD85]{Ast-K3}
A.~Beauville, J.-P. Bourguignon, and M.~Demazure (eds.), \emph{{G}\'eom\'etrie
  des {S}urfaces {K}3: {M}odules et {P}\'eriodes}, Ast\'erisque, vol. 126, Soc.
  Math. de France, 1985.

\bibitem[BLR90]{NM}
S.~Bosch, W.~L\"utkebohmert, and M.~Raynaud, \emph{{N}\'eron {M}odels},
  Springer-Verlag, 1990.

\bibitem[Del72]{D-K3}
P.~Deligne, \emph{{L}a {C}onjecture de {W}eil pour les {S}urfaces {K}3},
  Invent. Math. \textbf{15} (1972), 206--226.

\bibitem[Del81a]{Del-CanCoord}
\bysame, \emph{{C}ristaux {O}rdinaires et {C}oordonn\'ees {C}anoniques}, in
  Surfaces Alg\'ebriques, Lecure Notes in Mathematics, vol. 868,
  Springer-Verlag, 1981, pp.~80--127.

\bibitem[Del81b]{Del-K3}
\bysame, \emph{{R}el\`evement des {S}urfaces {K}3 en {C}aract\'etistique
  {N}ulle}, in Surfaces Alg\'ebriques, Lecture Notes in Mathematics, vol. 868,
  Springer-Verlag, 1981, pp.~58--79.

\bibitem[DM68]{D-M}
P.~Deligne and D.~Mumford, \emph{{T}he {I}rreducibility of the {S}pace of
  {C}urves of {G}iven {G}enus}, Publ. I.H.E.S. \textbf{36} (1968), 75--109.

\bibitem[Fri84]{Fri-Torelli}
R.~Friedman, \emph{{A} {N}ew {P}roof of the {G}lobal {T}orelli {T}heorem for
  {K}3 {S}urfaces}, Ann. Math. \textbf{120} (1984), 237--269.

\bibitem[GD67]{EGA}
A.~Grothendieck and J.~Dieudonn\'e, \emph{{\'E}l\'ements de {G}\'eom\'etrie
  {A}lg\'ebrique}, vol. 4,8,11,17,20,24 and 32, Publ. Math. de l'I.H.E.S.,
  Paris, 1960-1967.

\bibitem[Gro62]{FGA}
A.~Grothendieck, \emph{{F}ondaments de la {G}\'eom\'etrie {A}lg\'ebrique},
  Secr\'etariat Math\'ematique, Paris, 1962.

\bibitem[Har77]{HAG}
R.~Hartshorne, \emph{{A}lgebraic {G}eometry}, Graduate Texts in Mathematics,
  vol.~52, Springer-Verlag, New-York, 1977.

\bibitem[KM97]{K-M}
S.~Keel and S.~Mori, \emph{{Q}uotients by {G}roupoids}, Ann. of Math.
  \textbf{145} (1997), 193--213.

\bibitem[Knu71]{Knu-AS}
D.~Knutson, \emph{{A}lgebraic {S}paces}, Lecture Notes in Mathematics, vol.
  203, Springer-Verlag, 1971.

\bibitem[LMB00]{L-MB}
G.~Laumon and L.~Moret-Bailly, \emph{{C}hamps {A}lg\'ebriques}, Ergebnisse der
  Mathematik und ihrer Grenzgebiete, no.~39, Springer, 2000.

\bibitem[LP81]{L-P}
E.~Looijenga and C.~Peters, \emph{{T}orelli {T}heorems for {K}\"ahler {K}3
  {S}urfaces}, Comp. Math. \textbf{42} (1981), 145--186.

\bibitem[Mat58]{Mat}
T.~Matsusaka, \emph{{P}olarized {V}arieties, {F}ields of {M}oduli and
  {G}eneralized {K}ummer {V}arieties of {P}olarized {A}belian {V}arieties},
  Amer. Jour. of Math. \textbf{80} (1958), 45--82.

\bibitem[Mil80]{Mil-EC}
J.~S. Milne, \emph{{E}tale {C}ohomology}, Princeton University Press,
  Princeton, 1980.

\bibitem[MM64]{M-M}
T.~Matsusaka and D.~Mumford, \emph{{T}wo {F}undamental {T}heorems on
  {D}eformations of {P}olarized {V}arieties}, Amer. Jour. of Math. \textbf{86}
  (1964), 668--684.

\bibitem[Mum65]{Mum-GIT}
D.~Mumford, \emph{{G}eometric {I}nvariant {T}heory}, Ergebnisse der Mathematik
  und ihrer Grenzgebiete, no.~34, Springer, 1965.

\bibitem[Mum74]{Mum-AV}
\bysame, \emph{{A}belian {V}arieties}, Oxford University Press, New-York, 1974.

\bibitem[Nik80]{Nik-QLat}
V.~V. Nikulin, \emph{{I}ntegral {S}ymmetric {B}ilinear {F}orms and {S}ome of
  {T}heir {A}pplications}, Math. USSR Izv. \textbf{14} (1980), 103--166.

\bibitem[Ogu79]{O-SK3}
A.~Ogus, \emph{{S}upersingular {K}3 {S}urfaces}, Ast\'erisque \textbf{64}
  (1979), 3--86.

\bibitem[Ols04]{Olsson-K3}
M.~Olsson, \emph{{S}emi-{S}table {D}egenerations and {P}eriod {S}paces for
  {P}olarized {K}3 {S}urfaces}, Duke Math. J. \textbf{125} (2004), no.~1,
  121--203.

\bibitem[Oor62]{O-SchDeP}
F.~Oort, \emph{{S}ur le sch\'ema de {P}icard}, Bull. Soc. Math. Fr. \textbf{90}
  (1962), 1--14.

\bibitem[PSS72]{PSh-Sh}
I.~I. Piatetskij-Shapiro and I.~R. Shafarevich, \emph{{A} {T}orelli {T}heorem
  for {A}lgebraic {S}urfaces of {T}ype {K}3}, Math. USSR, Izv. \textbf{5}
  (1972), no.~3, 547--588.

\bibitem[Riz05]{JR-Thesis}
J.~Rizov, \emph{{M}oduli of {K}3 {S}urfaces and {A}belian {V}ariaties}, Ph.D.
  thesis, University of Utrecht, 2005.

\bibitem[RS76]{RSh-IMorK3}
A.~N. Rudakov and I.~R. Shafarevich, \emph{{I}nseperable {M}orphisms of
  {A}lgebraic {S}urfaces}, Math. USSR Izv. \textbf{10} (1976), no.~6,
  1205--1237.

\bibitem[SD74]{Saint-D}
B.~Saint-Donat, \emph{{P}rojective {M}odels of {K}-3 {S}urfaces}, Amer. Jour.
  of Math. \textbf{96} (1974), no.~4, 602--639.

\bibitem[Ser73]{S-CA}
J-P. Serre, \emph{{A} {C}ourse in {A}rithmetic}, Graduate Texts in Mathematics,
  vol.~7, Springer-Verlag, New York, 1973.

\bibitem[Shi79]{Shi-SSK3}
T.~Shioda, \emph{{S}upersingular {K}3 {S}urfaces}, in Algebraic Geometry, Proc.
  Summer Meeting, Univ. of Copenhagen, Lecture Notes in Mathematics,
  Springer-Verlag, 1979, pp.~564--591.

\bibitem[SI77]{Sh-I}
T.~Shioda and H.~Inose, \emph{{O}n {S}ingular {K}3 {S}urfaces}, in Complex
  Analysis and Algebraic Geometry, Iwanami Shoten/Cambridge University Press,
  1977, pp.~119--136.

\bibitem[Tat95]{JT-Cycles}
J.~Tate, \emph{{C}onjectures on {A}lgebraic {C}yclec in $l$-adic {C}ohomology},
  in Motives, S. Kleiman, U. Jansen and J.-P. Serre eds., vol.~55, Proc. of
  Symp. in Pure Math., no.~1, AMS, 1995, pp.~71--83.

\bibitem[Vie95]{Vie-M}
E.~Viehweg, \emph{{Q}uasi-projective {M}oduli for {P}olarized {M}anifolds},
  Ergebnisse der Mathematik und ihrer Grenzgebiete, vol.~30, Springer, 1995.

\end{thebibliography}
\providecommand{\bysame}{\leavevmode\hbox to3em{\hrulefill}\thinspace}
\providecommand{\MR}{\relax\ifhmode\unskip\space\fi MR }
\providecommand{\MRhref}[2]{%
  \href{http://www.ams.org/mathscinet-getitem?mr=#1}{#2}
}
\providecommand{\href}[2]{#2}

\end{document}